\documentclass[a4paper,12pt]{amsart}
\usepackage{ifthen,verbatim}
\usepackage{mathrsfs}
\usepackage{amsmath,amssymb,amsthm}

\usepackage{epsfig}
\usepackage{epstopdf}
\usepackage{pstricks}
\usepackage{pst-node}
\usepackage{pst-tree}
\usepackage{pst-plot}
\usepackage{pst-text}
\usepackage{multido}
\usepackage{subfigure}

\definecolor{hylightgray}{rgb}{0.9,0.9,0.9}

\numberwithin{equation}{section}
\nonstopmode
\theoremstyle{plain}
\newtheorem{cor}[equation]{Corollary}

\newtheorem{lem}[equation]{Lemma}
\newtheorem{prop}[equation]{Proposition}

\newtheorem{thm}[equation]{Theorem}
\theoremstyle{definition}
\newtheorem{rem}[equation]{Remark}

\newtheorem{nonsec}[equation]{}

\newtheorem{defn}[equation]{Definition}

\newenvironment{pf}[1][]{%
 \vskip 3mm
 \noindent
 \ifthenelse{\equal{#1}{}}%
  {{\slshape Proof. }}%
  {{\slshape #1.} }%
 }%
{\qed\bigskip}

\newcounter{alphabet}

\newenvironment{Thm}[1][]{\refstepcounter{alphabet}%
\bigskip%
\noindent%
{\bf Theorem \Alph{alphabet}}%
\ifthenelse{\equal{#1}{}}{}{ (#1)}%
{\bf .}
\itshape}{\vskip 8pt}

\newcommand{\be}{\begin{eqnarray}}
\newcommand{\ee}{\end{eqnarray}}
\newcommand{\ba}{\begin{array}}
\newcommand{\ea}{\end{array}}
\newcommand{\ben}{\begin{eqnarray*}}
\newcommand{\een}{\end{eqnarray*}}

\newcommand{\dist}{{\operatorname{dist}}}
\newcommand{\diam}{{\operatorname{diam}}}
\newcommand{\card}{{\operatorname{card}\,}}
\newcommand{\C}{{\mathbb C}}
\newcommand{\R}{{\mathbb R}}

\newcommand{\bB}{\mathbb{B}}

\newcommand{{\tth}}{\mathrm{th}}

\newcommand{\arth}{\,\mathrm{artanh}\,}
\newcommand{\capa}{\mathrm{cap}\,}
\newcommand{\D}{{\mathbb D}}
\newcommand{\ep}{\varepsilon}
\newcommand{\inv}{^{-1}}
\font\fFt=eusm10 %scaled 1200
\font\fFa=eusm7 %scaled 1200
\font\fFp=eusm5 %scaled 1200
\def\K{\mathchoice
{\hbox{\,\fFt K}}
{\hbox{\,\fFt K}}
{\hbox{\,\fFa K}}
{\hbox{\,\fFp K}}}
\newcommand{\Lam}{\Delta}

\newcommand{\mb}{\mathbb}

\newcommand{\pa}{\partial}

\newcommand{\sphere}{{\overline{\mathbb C}}}
\newcommand{\sig}{\sigma}
\newcommand{\Sig}{\Sigma}
\newcommand{\sq}{\sqrt}

\newcommand{\vep}{\varepsilon}
\newcommand{\vphi}{\varphi}

\newcommand {\Sn} {{\overline{\mathbb R}^n}}% the extended Euclidean n-space
\newcommand {\Stwo} {{\overline{\mathbb R}^2}}% the extended Euclidean n-space
\newcommand {\Rn} {{\mathbb R}^n}
\newcommand {\Rtwo} {{\mathbb R}^2}
\newcommand {\Hn} {{\mathbb H}^n}
\newcommand {\Bn} {{\mathbb B}^n}
\newcommand {\B} {{\mathbb B}}
\newcommand {\Btwo} {{{\mathbb B}^2}}
\newcommand {\Bbn} {{\overline{\mathbb B}^n}}
\newcommand {\Bber} {{\overline{B}}}
\newcommand{\mM}{\mathsf{M}}

\renewcommand {\mod} {{\rm mod}\,}

\newcommand{\sj}{{\,\hat j}}
\newcommand{\sd}{{\hat d}}

\definecolor{hylightgray}{rgb}{0.9,0.9,0.9}

\newcounter{minutes}
\divide\time by 60
\newcounter{hours}
\multiply\time by 60
\addtocounter{minutes}{-\time}
\begin{document}
\title[Conformally invariant complete metrics]{
Conformally invariant complete metrics
}
%OLD: Conformally invariant complete metrics
% OLD: Conformally invariant metrics and  domains with perfect boundaries
%Older title: On M\"obius invariant metrics and quasiregular maps into domains with uniformly perfect boundaries

\def\thefootnote{}
\footnotetext{
\texttt{\tiny File:~\jobname .tex,
          printed: \number\year-\number\month-\number\day,
          \thehours.\ifnum\theminutes<10{0}\fi\theminutes}
}
\makeatletter\def\thefootnote{\@arabic\c@footnote}\makeatother

\author[T. Sugawa]{Toshiyuki Sugawa}
\address{Graduate School of Information Sciences,
Tohoku University, Aoba-ku, Sendai 980-8579, Japan}
\email{sugawa@math.is.tohoku.ac.jp}
\author[M. Vuorinen]{Matti Vuorinen}
\address{Department of Mathematics and Statistics, University of Turku, FI-20014 Turku, Finland}
\email{vuorinen@utu.fi}
\author[T. Zhang]{Tanran Zhang}
\address{Department of Mathematics, Soochow University, No.1 Shizi Street, Suzhou 215006, China}
\email{trzhang@suda.edu.cn}


\begin{abstract}
For a domain $G$ in the one-point compactification
$\overline{\mathbb{R}}^n = \R^n \cup \{ \infty\}$ of $\R^n, n \ge 2$,
we characterize the completeness of the modulus metric
$\mu_G$ in terms of a potential-theoretic thickness condition of $\partial G\,,$ Martio's $M$-condition \cite{Mar75}. Next, we prove that $\partial G$ is
uniformly perfect if and only if $\mu_G$ admits a minorant in terms of a M\"obius
invariant metric. Several applications to quasiconformal maps are given.
\end{abstract}
\thanks{
The authors were supported in part by JSPS KAKENHI Grant Number JP17H02847
and NSF of the Higher Education Institutions of Jiangsu Province, China, Grant Number 17KJB110015, and NSFC Grant Number 12001391.
}
\maketitle

%%%%%%%%%%%%%%%%%%%%%%%%%%%%%%%%%%%%%%%
%%%%%%%%%%%%%%%%%%%%%%%%%%%%%%%%%%%%%%%
%%%%%%%%%%%%%%%%%%%%%%%%%%%%%%%%%%%%%%%
\section{Introduction}
%%%%%%%%%%%%%%%%%%%%%%%%%%%%%%%%%%%%%%%
%%%%%%%%%%%%%%%%%%%%%%%%%%%%%%%%%%%%%%%
%%%%%%%%%%%%%%%%%%%%%%%%%%%%%%%%%%%%%%%

Conformal invariance is one of the key notions in the geometric  theory of conformal
and quasiconformal maps both in the plane $\R^2 = \C$ and in the Euclidean space
$\R^n, n\ge3\,.$ Most clearly this is visible in the study of metrics: the uniformization theorem
\cite{bm} and the hyperbolic (Poincar\'e) metric of the
unit disk in $\C$ provide a way to define
the hyperbolic metric in any plane domain $G$ with $\card(\C\setminus G)\ge 2.$
This method fails for  $n\ge 3 $
because by Liouville's theorem   \cite{gmp, res} conformal maps in dimensions $n\ge 3 $
are M\"obius transformations. A widely studied natural question is whether some other
methods would work and whether there are counterparts of the hyperbolic metric in subdomains $G$ of
$\R^n$ and what sort of invariance or quasi-invariance properties, if any, such metrics
might have in higher dimensions $n\ge 3.$ From the vast literature we mention
A.~F.~Beardon \cite{be,be2}, J.~Ferrand \cite{LF,fe,fe2,Fer97,fmv}, F.~W.~Gehring \cite{gp,go,gh}, D.A. Herron \cite{bh,h,him,hj,hmm},
M.~Vuorinen \cite{vu85,vubook,hkv}. The recent extensive  research on metrics in geometric
function theory has many faces: two examples are the  monograph \cite{JP} of M.~Jarnicki  and P.~Pflug which provides
an encyclopedic treatise  on invariant metrics of complex manifolds and the monograph of
A.~Papadopoulos  which lists twelve metrics recurrent in geometric function theory \cite[pp. 42-48]{papa} .

Our main goal is to study one of these metrics, {\it the modulus metric} of a domain
$G\subset \overline{\R}^n = \R^n \cup \{\infty\}, n\ge 2\,,$
denoted by $\mu_G(x,y), x,y \in G\,,$ see Sections 3 and 4 for definitions.
In the special case of the unit ball, the modulus metric $\mu_{\B^n}(x,y)$ has an explicit formula
in terms of the hyperbolic metric of the unit ball $\B^n$;
the case of  $\mu_{\B^2}(x,y)$ was studied already by H.~Gr\"otzsch \cite[p.72]{a}.
The conformal invariant $\mu_G(x,y)$ has found numerous applications \cite{vubook,hkv},
but still many fundamental questions remain open.
Very recently a problem due to J.~Ferrand \cite{fmv}, \cite[pp.294-295]{hkv} was solved as follows.

\begin{Thm}[\cite{bpo,ps,z}] \label{thm:bpsz}
\black{A homeomorphism $f:G \to G'\,,$ where $G$ and $G'$ are domains in $\R^n, n\ge 2,$ is an isometry
between $(G,\mu_G)$ and $(G',\mu_{G'})$ if and only if $f$ is conformal.}
\end{Thm}

As pointed out above, $\mu_{\B^n}(x,y)$ is closely related to the hyperbolic metric of
$\B^n.$ We next study conditions on the domain $G$ under which $\mu_G$ defines an intrinsic
metric of $G$ having properties similar to the hyperbolic metric.
It turns out that the geometry of this metric significantly depends
on the ``potential theoretic thickness" of the boundary, measured in terms
of the conformal capacity.
As is well known, the conformal capacity is very closely connected
with the moduli of curve families  \cite[Thm 5.2.3, p. 164]{gmp}, \cite[Theorem 9.6, p. 152]{hkv}.

If the boundary $\partial G$ is polar, i.e. if it has null conformal capacity $\capa(\pa G) = 0$,
then $\mu_G \equiv 0$; otherwise $\mu_G$ is a conformally invariant metric.
Even if $\capa(\pa G)>0,$ the
modulus metric $\mu_G$ might not reflect the intrinsic geometry of $G$ very precisely.
For instance,  a polar compact set $N \subset G$ is invisible for the modulus metric in the sense
that if $\capa N=0,$ then $\mu_G(x,y)=\mu_{G\setminus N}(x,y)$ for
$x,y\in G\setminus N.$ Therefore, it is meaningful to look for a condition on
$G$ so as to guarantee that $\mu_G$ is a complete metric.
We remark that a similar problem for the Kobayashi metric on domains in $\C^n$
is rather difficult (see, e.g., \cite{gau,pf}).

In connection with
this completeness property, we recall another notion on metric spaces.
A metric space $(X,m)$ is called {\it proper} \cite{brh} if the closed metric ball
$\{x\in X: m(x,a)\le r\}$ is compact whenever $a\in X$ and $r>0.$
This is equivalent to say that the open metric ball $\{x\in X: m(x,a)<r\}$ is
relatively compact for $a\in X$ and $r>0.$
Note that a proper metric space is locally compact and complete.
However, the converse is not true in general.
(Consider, e.g., $(X,m/(1+m))$ for a locally compact but non-compact
complete metric space $(X,m)$ such as $\R^n$ with the Euclidean metric.)

Our first result characterizes domains $G$ for which the metrics $\mu_G$ are  complete.
\begin{thm}\label{thm:muproper}
Let $G$ be a domain in $\Sn$ with $\partial G\ne\emptyset.$
Then the following conditions are equivalent:
\begin{enumerate}
\item[(i)]
$(G,\mu_G)$ is a proper metric space.
\item[(ii)]
$(G,\mu_G)$ is a complete metric space.
\item[(iii)]
$G$ is an M-domain.
That is to say, each boundary point $x$ of $G$ satisfies the M-condition.
%That is to say, the boundary $\pa G$ of $G$ does not satisfy
%the continuum criterion for every boundary point $x$ of $G$:
%$\mM(x,\Sn\setminus G)=\infty.$
\end{enumerate}
\end{thm}

The M-condition for $x\in\pa G$ was introduced by O.~Martio \cite{Mar75}%
\footnote{{\black The M-condition $\mM(x,\Sn\setminus G)=\infty$} was denoted by $M_x=\infty$ in Martio's paper \cite{Mar75}.}
in his study of potential theoretic regularity of the domain. {\black If this condition
holds for all  $x\in\pa G,$} the complement $\Sn\setminus G$
of $G$ is ``thick enough" at every point of $\partial G\,$
% i.e. it does not satisfy the continuum criterion at $x$
 \cite{Mar75}, \cite{MarSar78}.
%$\mM(x,\Sn\setminus G)=\infty.$
See Section 3 for definitions of those concepts and related properties.

Our second result refines further the case when $\mu_G$ is complete.
We assume now that the boundary
of a domain is {uniformly perfect} in the sense of Ch.~Pommerenke \cite{Pom79, Pom84}
{\black  --- in this case
the M-condition is valid, see Corollary \ref{UPandM}.}
This notion was  introduced by A.~F.~Beardon and Ch.~Pommerenke \cite{bp} for unbounded
closed sets in $\C$, but about the same time an
equivalent concept was studied by P.~Tukia and J.~V\"ais\"al\"a \cite{tv}
under the name ``homogeneously dense sets" in the setting of general metric spaces.
By definition, a compact set $E$ in
$\Sn$ with $\card(E)\ge 2$ is called {\it uniformly perfect}
if there exists a constant $c\in(0,1)$ such that $E$ meets the closed annulus
$cr\le |x-a|\le r$ whenever $a\in E\setminus\{\infty\}$ and $r\in(0,\diam(E)),$
where $\diam(E)$ denotes the Euclidean diameter of $E$ and set $\diam(E)=+\infty$
when $\infty\in E.$  In the planar case when $G\subset\Rtwo=\C,$
A.~F.~Beardon and Ch.~Pommerenke \cite{bp} gave another characterization
in terms of the hyperbolic and quasihyperbolic metrics $h_G(x,y)$ and $k_G(x,y)\,,$
resp. (see Section 2), and proved that
$\pa G$ is uniformly perfect if and only if there is a constant $b>0$ such that
\[h_G(x,y)\ge b k_G(x,y)\quad {\textrm{ for}\, all\,\, } x,y\in G\,.\]
Here we give an alternative characterization of uniform perfectness of $\pa G$ in terms of
intrinsic metrics which is valid in higher dimensions as well and, moreover, is applicable
to subsets of the M\"obius space. This characterization  requires that the modulus metric
be minorized by a M\"obius invariant metric $\delta_G\,,$ defined in terms of the absolute ratio
\ref{deltaDef} for all domains $G \subset \overline{\R}^n$ with ${\textrm{ card}} (\partial G)\ge 2\,.$
This metric was first introduced in \cite[pp.115-116]{vubook}  and, later on, studied by
P.~Seittenranta in his PhD thesis \cite{Seitten99} where he also suggested the name ``M\"obius
metric".

\begin{thm} \label{thm:muminor}
Let $G \subset \Sn$ be a domain with $\card(\pa G)\geq 2.$
Then $\pa G$ is uniformly perfect if and only if
there exists a constant $b>0$ such that for all $x,y\in G$ the inequality
\begin{equation}\label{eq:mu-delta}
\mu_G (x,y) \geq b \,\delta_G (x,y)
\end{equation}
holds, where $\mu_G$ is the modulus metric and $\delta_G$ is the M\"obius metric.
\end{thm}

For a proper subdomain $G$ of $\Rn,$  the lower bound \ref{eq:mu-delta}
can be expressed in terms of a similarity invariant metric,
the {\it distance-ratio metric} of $G$ as follows. For $x,y \in G$ define
\begin{equation}\label{jdef}
j_G(x, y) = \log\left(1+\frac{|x-y|}{\min\{d_G(x),d_G(y)\}}\right),
\end{equation}
which is a metric on $G,$ where $d_G(x)$ denotes the Euclidean distance from $x$
to the boundary $\pa G$ \cite[Lemma 4.6, p. 59]{hkv}.
When $G\subset\R^n,$ the above condition  \ref{eq:mu-delta} is equivalent to
the requirement that for some constant $b'>0$
\[\mu_G(x,y)\ge b' j_G(x,y)\]
for all $x,y \in G\,.$
Since $(G,\delta_G)$ is a proper metric space (see Lemma \ref{lem:proper} below),
we have the following result as a corollary of
Theorems \ref{thm:muproper} and \ref{thm:muminor}.

\begin{cor} \label{UPandM}
Let $G \subset \Sn$ be a domain with $\card(\pa G)\geq 2.$
If $\pa G$ is uniformly perfect, then $G$ is an M-domain.
\end{cor}

The converse is not true in general. A counterexample will be given in Section 3.

The proof of Theorem \ref{thm:muminor} is based, in part, on a potential theoretic
thickness characterization of uniform perfectness  \cite{vu}, \cite{JV96}.
Many authors have contributed to the research of
uniformly perfect sets and related thickness conditions,
see \cite{aw}, \cite{brc}, \cite[pp. 343-345]{gm}, \cite{gsv},  \cite{KL}, \cite{L}, \cite{marm},
\cite{mm} and the survey of T. Sugawa \cite{sug} on uniform perfectness.

Uniform domains play an important role in geometric function theory.
See \cite{go} and the recent monograph \cite{gh} for details.
For convenience of the reader, we will provide a brief account on this notion in the next section.

\begin{thm} \label{thm:mumajor}
Suppose that $G \subset \Sn$ is a uniform domain.
Then there exist constants $d_1, d_2$
depending only on $n$ and the uniformity parameters such that
\be \label{eq:d}
\mu_G (x,y) \leq d_1\, \delta_G (x,y)+d_2\, \quad x,y\in G.
\ee
Conversely, suppose that a domain $G$ in $\Stwo$ with continuum as its boundary
satisfies \ref{eq:d}. Then $G$ is uniform.
\end{thm}

Note that the boundary of a domain $G$ in $\Stwo=\sphere$ is a continuum;
that is, a non-degenerate connected compact set,
if and only if $G$ is a simply connected hyperbolic domain.
It is known that such a domain $G$ is uniform precisely when
$G$ is a quasidisk,
that is to say, $G$ is the image of the unit disk $\bB^2$ under a quasiconformal
homeomorphism of $\sphere$ onto itself \cite{gh}.
Therefore, as a corollary, we have the following characterization of quasidisks.

\begin{cor} \label{mycor}
Let $G$ be a simply connected domain in the Riemann sphere $\sphere$
with $\card(\sphere\setminus G)\ge 2.$
Then $G$ is a quasidisk if and only if there are positive constants
$d_1$ and $d_2$ such that the inequality
$$\mu_G (z,w) \leq d_1\, \delta_G (z,w)+d_2$$
holds for all $z,w\in G.$
\end{cor}

In this corollary, we may replace the modulus metric $\mu_G$ by
Ferrand's modulus metric $\lambda_G^{-1}$ (see Lemma \ref{lem:mu-lambda} below).
We remark that for  $G\subset \C$ the above condition is also equivalent to the condition
 \[\mu_G(z,w)\le d_1'\, j_G(z,w)+d_2' \quad {\textrm{ for} \,\,} z,w\in G\,.\]
As we will see later, the constant $d_2$ in Corollary \ref{mycor} cannot be dropped.
We expect that the converse would be true for all dimensions $n\ge2$
under a weaker assumption on the boundary such as uniform perfectness of the boundary.
{\black
These observations lead to the following problem.

\begin{nonsec}{\bf Open problem.}
{\rm
Let $n\ge 2.$
Find a geometric condition $(*)$ on the boundaries of domains $G$ in $\Sn$ with the
following property: If a domain $G$ in $\Sn$ satisfies the condition $(*)$ and
the inequality \ref{eq:d} for some constants $d_1>0$ and $d_2>0,$ then $G$ is uniform.
}
\end{nonsec}
}

Finally, we consider the hyperbolic metric $h_G$  {\black{and the Ferrand metric $\sigma_G,$
see \ref{eq:w},}} in planar domains $G.$
It is well known \cite{bp}  that %\cite{gh}, \cite[Ch 15]{KL},
if $\partial G$ is uniformly perfect, then the distances in the $h_G$ metric are comparable to those in the quasihyperbolic metric $k_G\,.$
{\black Furthermore,} this comparison property fails to hold if the domain $G$ has isolated boundary points.
Indeed, the following asymptotic formulae hold.

\begin{lem}\label{lem:asymp}
Let $G$ be a hyperbolic domain in $\sphere$
and suppose that $G$ has an isolated boundary point $a$ with $a\ne\infty.$
Then, for a fixed $z_0\in G,$ as $z\to a$
\begin{equation}\label{eq:asymp1}
\sigma_G(z,z_0)=\log\frac1{|z-a|}+O(1)
\quad\textrm{and}\quad
\delta_G(z,z_0)=\log\frac1{|z-a|}+O(1),
\end{equation}
while
\begin{equation}\label{eq:asymp2}
h_G(z,z_0)=\log\log\frac1{|z-a|}+O(1).
\end{equation}
\end{lem}

It is a challenging task, studied in \cite{SVZ19} and \cite{SZ17},
to give concrete bounds for the $h_G$ distances in domains $G$ whose boundary
consists only of isolated points.
Since $\log(1+x)$ is a subadditive function on $0\le x<+\infty,$
we can easily see that $\log(1+m(x,y))$ is a distance function on $X$ whenever
$m(x,y)$ is a distance function on $X$ \cite[7.42(1)]{avv}.
In view of the above behaviour of the hyperbolic distance around isolated boundary points,
we are led to the introduction of the {\it logarithmic M\"obius metric} $\Lam_G(x,y)$
and the {\it logarithmic Ferrand metric} $\Sig_G(x,y)$ for a domain
$G \subset \Sn $ with $\card(\Sn \setminus G)\geq 2$ as follows:
\begin{eqnarray} \label{logMm}
\Lam_G(x,y)&= \log(1+ \delta_G(x,y))\,,\quad x,y \in G\,, \\
\Sig_G(x,y)&= \log(1+ \sig_G(x,y))\,,\quad x,y \in G\,.
\label{logsig}
\end{eqnarray}
Because $ \delta_G$ and $\sig_G$ are M\"obius invariant, $\Lam_G$ and $\Sig_G$ are M\"obius invariant metrics, too.
We also have $\Lam_G(x,y)\le \Sig_G(x,y)$ (see Lemma \ref{lem:comparedelta} below).
When the complement of $G$ in $\sphere$ is a finite set, the hyperbolic distance
$h_G$ is majorized by $\Lam_G.$
However, $h_G$ is never minorized by it for any domain with a puncture;
namely, with an isolated boundary point.
In fact, we prove a slightly stronger result.

\begin{thm} \label{thm:hypmajor}
Let $A$ be a finite set in $\sphere$ with $\card(A)\ge 3$
and let $G=\sphere\setminus A.$
Then there exists a positive constant $c=c(A)$ such that for all $z,w \in G$,
$$h_{G}(z,w) \leq c\, \Lam_{G} (z,w)=c \log(1+ \delta_G(z,w))\,.$$
On the other hand, for an arbitrary hyperbolic domain $G$ in $\sphere$
with a puncture,
there is no non-decreasing function $\Phi:[0,+\infty)\to[0,+\infty)$
with $\Phi(t)>0$ for $t>0$ such that for all $z,w\in G,$
$$
\Phi(\delta_G(z,w))\le h_G(z,w).
$$
\end{thm}

All the results here will be proved in the subsequent sections.
More precisely, this paper is organized as follows.
Section 2 is devoted to definitions and basic properties
of the metrics involved, with the exception of the modulus metric, which will be defined
in Section 4.
In Section 3, we recall the notion of the (conformal) modulus of
a curve family and its fundamental properties.
We also introduce the notion of M-domains defined in terms of the
continuum criterion of Martio \cite{Mar75}.
The modulus metric is defined and related results are established in Section 4.
We give some applications of the above results to quasiconformal or quasiregular
mappings in Section 5. Theorem \ref{thm:hypmajor} is proved in the last section.
Two open problems are pointed out, namely \ref{HdimQ} and \ref{HypQ}.

\bigskip

%%%%%%%%%%%%%%%%%%%%%%%%%%%%%
%%%%%%%%%%%%%%%%%%%%%%%%%%%%%
%%%%%%%%%%%%%%%%%%%%%%%%%%%%%

\section{Preliminary notation and results}
%%%%%%%%%%%%%%%%%%%%%%%%%%%%%
%%%%%%%%%%%%%%%%%%%%%%%%%%%%%
%%%%%%%%%%%%%%%%%%%%%%%%%%%%%

We follow standard notation. See e.g. \cite{be}, \cite{va} for more details.
We write
\begin{eqnarray*}
B^{n} (x,r)&=\{z \in \Rn: |z-x| < r\}, \\
\Bber^{n} (x,r) &=\{z \in \Rn: |z-x| \le r\}, \\
S^{n-1} (x,r)&=\{z \in \Rn: |z-x| = r\},
\end{eqnarray*}
for balls and spheres, respectively, and
$$\mb{B}^{n}=B^{n} (0,1),\quad \mb{H}^{n} = \{(x_1,\ldots, x_n) \in \Rn: x_n > 0 \}.$$

First we recall the definition of the {\it chordal (spherical)
distance} $q(x,y)$ on $\Sn$:

\be\label{def:chordal}\left\{
\begin{array}{ll}
\displaystyle
q(x,y)=\frac{|x-y|}{\sqrt{1+|x|^2}\sqrt{1+|y|^2}}\,, &\quad ~x,y \neq \infty\,, \\[5mm]
\displaystyle q(x,\infty)=q(\infty,x)=\frac{1}{\sqrt{1+|x|^2}}\,, &\quad ~x \neq \infty\,.
\end{array}
\right.
\ee

For distinct points $a, b, c, d \in \Sn$, the {\it absolute (cross) ratio} is defined by
$$|a, b, c, d|=\frac{q(a,c)q(b,d)}{q(a,b) q(c,d)}\,.$$
When none of the points is $\infty,$ we see that
$$|a, b, c, d|=\frac{|a-c||b-d|}{|a-b||c-d|}\,.$$

\medskip

\begin{nonsec}{\bf Hyperbolic metric.} \label{defHyp}{\rm
The hyperbolic metrics $2|dx|/(1-|x|^2)$ on $\Bn$ and $|dx|/x_n$ on $\Hn$
induce the hyperbolic distances $h_{\Bn}(x,y)$ and $h_{\Hn}(x,y)$ respectively.
When $n=2,$ any domain $G$ of $\overline\R^2=\sphere$ with $\card(\partial G)\ge 3$
is known to have a holomorphic universal covering projection $p$ of the unit disk $\mb{B}^2$ onto $G.$
Thus the hyperbolic distance $h_G$ of $G$ can be defined by
$$
h_G(z_1,z_2)=\min_{\zeta_1\in p\inv(z_1), \zeta_2\in p\inv(z_2)} h_{\mb{B}^2}(\zeta_1,\zeta_2)
=\inf_{\gamma\in\Gamma} \int_\gamma \rho_G(z)|dz|,
$$
where $\Gamma$ is the set of all rectifiable curves joining $z_1$ and $z_2$ in $G$
and $\rho_G(z)$ denotes the hyperbolic density determined by the relation
$2/(1-|\zeta|^2)=\rho(p(\zeta))|p'(\zeta)|,~ \zeta\in\mb{B}^2$
(see \cite{bm}, \cite{KL} for details).}
\end{nonsec}

\begin{nonsec}{\bf Quasihyperbolic metric.} \label{defQHyp}{\rm
For higher dimensions, however, we cannot define hyperbolic metric for general domains.
Quasihyperbolic metrics were introduced by F.W. Gehring and B. Palka \cite{gp} as a substitute for it.
For a domain $G\subsetneq {\mb{R}^n}$, the {\it quasihyperbolic metric} ${k}_G$ is defined by
$${k}_{G}(x,y)=\inf_{\gamma \in \Gamma} \int_{\gamma}\frac{|dt|}{d_G(t)}, \quad x, y \in G,$$
where $\Gamma$ is the family of all rectifiable curves in $G$ joining $x$ and $y$.
Note here that the inequality
$$j_G(x,y)\le k_G(x,y)$$
holds for an arbitrary $G\subsetneq {\mb{R}^n}$ and all
$x,y \in G$ \cite[Lemma 2.1]{gp}.}
\end{nonsec}

\begin{nonsec}\label{unifDom}{\bf Uniform domains.} {\rm
A proper subdomain $G$ of $\Rn$ is called {\it uniform}
if there exist positive constants $a$ and $b$ with
the following property \cite{MarSar78, go}:
for every pair of points $x_1, x_2\in G,$ there is a rectifiable curve $\gamma$
joining $x_1$ and $x_2$ in $G$ in such a way that
$\ell(\gamma)\le a|x_1-x_2|$ and that $\min\{\ell(\gamma_1), \ell(\gamma_2)\}\le b\, d_G(x)$
for each $x\in\gamma,$ where $\gamma_j$ is the part of $\gamma$ between $x_j$ and $x$
for each $j=1,2,$ $\ell(\gamma)$ denotes the length of the curve $\gamma$ and
$d_G(x)$ is the Euclidean distance to the boundary of $G$ from $x.$
The class of uniform domains can also be defined in terms of a comparison
inequality between two metrics \cite{go,vubook}
\footnote{In \cite{go}, the condition \ref{defUnif} was given in the slightly different form
$k_G(x,y)\le a\, j_G(x,y)+b$ for some constants $a,b.$
We easily see that we can take $b=0$ by letting $a$ be larger if necessary.
See \cite[2.50 (2)]{vu85}.}
a subdomain $G$ of $\Rn$ with non-empty boundary
is uniform if and only if there exists a constant $c\ge 1$ such that
\begin{equation} \label{defUnif}
k_G(x,y) \le c\, j_G(x,y)
\end{equation}
for all $x,y \in G$, where $k_G$ and $j_G$ are the quasihyperbolic
and distance-ratio metrics, respectively. Note that $j_G(x,y)\le k_G(x,y)$
holds for every domain $G$ and all $x,y \in G$ by \cite[Lemma 2.1]{gp}.}
\end{nonsec}

\begin{nonsec}{\bf Ferrand's metric.} \label{defFerrand}{\rm
Since the definition of the quasihyperbolic metric relies on the Euclidean metric, it is not defined
for all subdomains of the M\"obius space and therefore
it is not  M\"obius invariant. To overcome this shortcoming{\black{,}}
J.~Ferrand \cite{fe}  modified the definition as follows.
For a subdomain $G$ of $\Sn$ with $\card (\pa G) \geq 2,$ define a density function
$$
w_G (x)=\sup_{a,b \in \pa G} \frac{|a-b|}{|x-a|\,|x-b|},\quad x \in G \setminus \{\infty\}\,,
$$
and the metric $\sig_G$ in $G$,
\begin{equation}\label{eq:w}
\sig_{G}(x,y)=\inf_{\gamma \in \Gamma}\int_{\gamma} w_G (t) |dt|,
\end{equation}
where $\Gamma$ is the family of all rectifiable curves in $G$ joining $x$ and $y$.
The following result is due to Ferrand \cite[p. 122]{fe} and $\sig_G(x,y)$
is now called the {\it Ferrand metric} \cite[Ch. 5]{hkv}.}
\end{nonsec}

\begin{lem}\label{lem:sigandkap}
Let $G \subset {\Sn}$ be a domain with $\card (\pa G) \geq 2$.
The Ferrand metric $\sig_{G}$ has the following properties.
\begin{enumerate}
\item
$\sig_{G}$ is a M\"obius invariant metric.
\item
When $ G $ is either $\Bn$ or $\Hn,$
 $\sig_{G}$  coincides with the hyperbolic metric $h_{G}\,.$
\item
${k}_G \leq \sig_G \leq 2 {k}_G$ for every domain $G \subsetneq \Rn$.
\end{enumerate}
\end{lem}

We remark that the metric $\sigma_G$ was recently studied by
D.~A.~Herron and P.~K.~Julian \cite{hj}.

\begin{nonsec}{\bf M\"obius metric.} \label{defMobius}{\rm
Let $G \subset \Sn$ be an open set with $\card (\pa G) \geq 2$.
The {\it M\"obius metric}  on $G$ is defined as follows  (\cite[pp.115-116]{vubook}, Seittenranta\cite{Seitten99}):
\begin{equation} \label{deltaDef}
\delta_G(x,y):=\log (1+ m_G (x,y))\,, \quad m_G (x,y):=\sup_{a,b \in \pa G} |a,x,b,y|\,.
\end{equation}
%We call it the M\"obius metric of Seittenranta.
%The following lemma is due to Vuorinen \cite[Lemma 8.39]{vubook}.
Note that the M\"obius metric $\delta_{G}$  coincides with the hyperbolic metric $h_{G}$
when $ G $ is either $\Bn$ or $\Hn$ \cite[Lemma 8.39]{vubook}. A metric very similar to
the M\"obius metric is the Apollonian metric of A. F. Beardon \cite{be2}.}
\end{nonsec}

\begin{nonsec}{\bf Chordal distance-ratio metric.} \label{defjhat}
{\rm
For a proper subdomain $G$ of $\Sn$ we define the {\it chordal (spherical) distance-ratio metric} by
$$
\sj_G(x, y) = \log\left(1+\frac{q(x,y)}{\min\{\sd_G(x), \sd_G(y)\}}\right),
$$
where
$$
\sd_G(x)=\inf_{a\in\partial G}q(x,a).
$$
The triangle inequality for this metric follows from \cite[Lemma 2.2]{Seitten99}.
}
\end{nonsec}

The following results are due to Seittenranta \cite{Seitten99}.

\begin{lem}%{ \cite[Theorems 3.3, 3.4, 3.6, 3.12]{Seitten99}}
\label{lem:comparedelta}
Let $G$ be an open subset of $\Sn$ with $\card (\pa G) \geq 2\,.$
Then $\delta_{G}$ is a M\"obius invariant metric and the following hold:
\begin{enumerate}
\item
$\delta_G \leq  \sig_G$ .
\item
$\delta_G \le 2\sj_G \,.$
\item
If $G \subsetneq \Rn\,,$ then $j_G \leq \delta_G \leq 2 j_G\,.$
\end{enumerate}
\end{lem}

\begin{pf}
The fact that $\delta_G$ satisfies the triangle inequality, assertions
(1) and (3) follow
 from Theorems 3.3, 3.4 and 3.12 in \cite{Seitten99}, respectively.
 In order to show assertion (2), we introduce the auxiliary metric
 $$
 j^*_G(x,y)=
 \log\left(1+\frac{q(x,y)}{\sd_G(x)}\right)
 +\log\left(1+\frac{q(x,y)}{\sd_G(y)}\right).
 $$
 Theorem 3.6 in \cite{Seitten99} means the inequality $\delta_G(x,y)\le
 j^*_G(x,y)$
 for $x,y\in G.$
 It is easy to verify
 the inequalities $\sj_G(x,y)\le j^*_G(x,y)\le 2\sj_G(x,y).$
 Thus assertion (2) now follows.
 \end{pf}
\hfill  $\Box$

As a consequence of the previous lemma, we have the following inequality,
which will be used in the proof of Theorem \ref{thm:muminor} later:
\begin{equation}\label{eq:j2j}
j_G(x,y)\le 2\sj_G(x,y),\quad x,y\in G\subsetneq\R^n.
\end{equation}
We note that there is no constant $c=c(n)>0$ depending only on $n$ such that
$j_G(x,y)\ge c\sj_G(x,y),~ x,y\in G,$ holds for all proper subdomains $G$
of $\R^n.$
The following result follows also from the previous lemma.

\begin{lem}\label{lem:proper}
The metric space $(G,\delta_G)$ is proper for $G\subset\Sn$ with $\card(\pa G)\ge 2.$
\end{lem}

\begin{pf}
By the M\"obius invariance, we may assume that $G\subset\R^n.$
Then $j_G\le \delta_G$ by Lemma \ref{lem:comparedelta} (3).
Therefore, it is enough to show that $(G,j_G)$ is proper in this case.
For $a\in G$ and $0<r,$ we have to show that the set $B=\{x\in G: j_G(x,a)<r\}$ is relatively compact.
It is enough to show that $B$ is bounded and $\dist(B, \pa G)>0.$
The inequality $\log(1+|x-a|/d_G(a))\le r$ holds for $x\in B$ and thus
$|x-a|\le d_G(a)(e^r-1),$ which proves that $B$ is bounded.
On the other hand, the inequality $\log(1+|x-a|/d_G(x))\le r$ holds for $x\in B.$
Note that $d_G(x)\ge d_G(a)/2$ if $|x-a|\le d_G(a)/2.$
For $x\in B$ with $|x-a|\ge d_G(a)/2,$ we thus have
$d_G(x)\ge |x-a|/(e^r-1)\ge d_G(a)/(e^r-1).$
Therefore, we have shown $\dist(B,\pa G)\ge \min\{d_G(a)/2, d_G(a)/(e^r-1)\}>0$
as required. \hfill $\Box$
\end{pf}

\begin{nonsec}{\bf M\"obius uniform domains.}{\rm
We now consider a M\"obius invariant characterisation of uniform domains.
As we saw above, uniform domains in $\R^n$ are characterised by the condition
\ref{defUnif} in terms of quasihyperbolic and distance-ratio metrics.
These two metrics are invariant under similarity transformations but unfortunately
not under M\"obius transformations. To overcome this lack of invariance we apply Ferrand's
M\"obius invariant metric $\sigma_G$ and the M\"obius metric $\delta_G\,.$}
\end{nonsec}

\begin{defn}[\cite{Seitten99}] \label{MobUnif}{\rm
We say that a domain  $G \subset \Sn\,$ with $\card(\Sn\setminus G)\ge 2$
is {\it M\"obius uniform},
if there exists a constant $c\ge1$ such that for all $x,y \in G$
$$
\sigma_G(x,y)\le c \,\delta_G(x,y)\,.
$$}
\end{defn}

Note that definition \ref{defUnif} only applies to subdomains of $\Rn$ whereas
Definition \ref{MobUnif} applies to subdomains of $ \Sn .$
Indeed, we have the following result.

\begin{prop}\label{prop:EqUnif}
Let $G \subset \Rn$ be a domain with $\card (\pa G) \geq 2$.
Then $G$ is M\"obius uniform if and only if it is uniform in the sense of \ref{defUnif}.
\end{prop}

\begin{pf} From Lemmas \ref{lem:sigandkap} and \ref{lem:comparedelta} it follows that if $G$ is
M\"obius uniform with a constant $c_1$, then it is uniform in the sense of \ref{defUnif}
with the constant $2 c_1\,.$
Conversely, from Lemmas \ref{lem:sigandkap} and \ref{lem:comparedelta} it follows that if $G$ is
 uniform in the sense of \ref{defUnif} with a constant $c_2$, then it is M\"obius uniform
with the the constant $2 c_2\,.$ \hfill $\Box$
\end{pf}

Therefore, we will use the shorter term ``uniform" below for both
uniform domains and M\"obius uniform domains
unless we want to emphasize which definition is intended.

We end this section with a proof of Lemma \ref{lem:asymp}.

\begin{pf}[Proof of Lemma \ref{lem:asymp}]
By assumption, there is a number $r>0$ such that the punctured disk $0<|z-a|<r$
is contained in $G.$
It is enough to prove the assertions for $a=0$ and $r=1.$
By assumption, we can find a finite boundary point $b$ of $G$ so that
$$
m_G(z,z_0)\ge |0,z,b,z_0|=\frac{|b||z-z_0|}{|z||b-z_0|}
\ge\frac{|b||z_0|}{2|z||b-z_0|}=:\frac{C}{|z|}
$$
for $z\in G$ with $0<|z|<|z_0|/2.$
Hence,
$$
\delta_G(z,z_0)=\log(1+m_G(z,z_0))\ge\log(1+C/|z|)=\log\frac1{|z|}+O(1)
$$
as $z\to 0.$
Next we estimate $w_G(z)$ from above for $0<|z|\le 1/4.$
For $b\in\partial G\setminus\{0\},$ we have
$|z-b|/|b|\le 1+|z|/|b|\le 1+|z|$ and
$|z-b|/|b|\ge 1-|z|/|b|\ge 1-|z|$ and thus
$$
%\frac{|a-b|}{|z-a||z-b|}=
\frac{16}5\le\frac1{|z|(1+|z|)}\le
\frac{|b|}{|z||z-b|}\le \frac1{|z|(1-|z|)}
=\frac1{|z|}+\frac1{1-|z|}
\le\frac1{|z|}+\frac43
$$
for $0<|z|\le 1/2.$
For $b_1,b_2\in\partial G\setminus\{0\},$ we have $|z-b_j|\ge |b_j|-|z|\ge 3|b_j|/4\ge 3/4$ and
$$
\frac{|b_1-b_2|}{|z-b_1||z-b_2|}
\le \frac{|z-b_2|+|z-b_1|}{|z-b_1||z-b_2|}
=\frac{1}{|z-b_1|}+\frac{1}{|z-b_2|}\le \frac{8}3
$$
as $z\to0.$
Hence, we obtain $w_G(z)\le 1/|z|+4/3$ for $0<|z|\le 1/4.$
For a given $z_0,$ we take a point $z_1\in G$ so that $|z_1|\le \min\{|z_0|,1/4\}.$
Then, for $0<|z|<|z_1|,$ we have
$$
\sigma_G(z,z_0)\le\sigma_G(z,z_1)+\sigma_G(z_1,z_0)
\le \int_\gamma \frac{|dt|}{|t|}+O(1)=\log\frac1{|z|}+O(1),
$$
where $\gamma$ is the curve going from $z_1$
to the point $(|z_1|/|z|)z$  along the circle $|t|=|z_1|$ and then going to $z$ radially.
Since $\delta_G(z,z_0)\le\sigma_G(z,z_0),$ \ref{eq:asymp1} follows.

Secondly, we prove \ref{eq:asymp2}.
For simplicity, we further assume that $1,\infty\in\partial G.$
(For the general case, we may use a suitable M\"obius transformation to reduce to this case.)
Then
$$
\D^*=\{z\in\C: 0<|z|<1\}\subset G\subset \C\setminus\{0,1\}
$$
and therefore
$$
\rho_{\D^*}(z)\ge\rho_G(z)\ge \rho_{\C\setminus\{0,1\}}(z)
$$
for $0<|z|<1.$
Since
$$
\rho_{\D^*}(z)=\frac1{|z|\log(1/|z|)}
\quad\textrm{and}\quad
\rho_{\C\setminus\{0,1\}}(z)=\frac1{|z|(C_0+\log(1/|z|))},
$$
where $C_0=1/\rho_{\C\setminus\{0,1\}}(-1)$ (see \cite{KL} for instance),
we have
$$
\rho_G(z)=\frac1{|z|\log(1/|z|)}+O\left(\frac1{|z|\log^2(1/|z|)}\right)
$$
as $z\to0.$
Noting the fact that the real function $1/[t\log^2t]$ is integrable over $(0,1/2],$
we obtain the required asymptotics \ref{eq:asymp2} as required. \hfill $\Box$
\end{pf}

\begin{rem}{\rm
As the above proof shows, \ref{eq:asymp1} is valid also in dimensions $n\ge 2.$}
\end{rem}

%%%%%%%%%%%%%%%%%%%%%%%%%%%%%%%%%%%%%%%%%%%%%%%%%%%%%%%%%%%%%%%%%%%
%%%%%%%%%%%%%%%%%%%%%%%%%%%%%%%%%%%%%%%%%%%%%%%%%%%%%%%%%%%%%%%%%%%
\section{Modulus and M-domains} \label{sec:M-domains}
%%%%%%%%%%%%%%%%%%%%%%%%%%%%%%%%%%%%%%%%%%%%%%%%%%%%%%%%%%%%%%%%%%%
%%%%%%%%%%%%%%%%%%%%%%%%%%%%%%%%%%%%%%%%%%%%%%%%%%%%%%%%%%%%%%%%%%%
We recapitulate some of the basic facts about moduli of curve families
and quasiconformal maps, following \cite{gmp,va}.
Let $\Gamma$ be a family of curves in $\Sn$.
We say that a non-negative Borel-measurable function $\rho: \Rn\to \R \cup \{+\infty\}$
is an admissible function for $\Gamma,$
if $\int_{\gamma} \rho ds \geq 1$ for each locally rectifiable curve $\gamma$ in $\Gamma$.
The (conformal) modulus of $\Gamma$ is
%% rmk: \mathfrak --> \mathcal
$$\mM(\Gamma)=\inf_{\rho \in \mathcal{F}(\Gamma)}\int_{\Rn}\rho^n dm,$$
where $\mathcal{F}(\Gamma)$ is the family of admissible functions for $\Gamma$
and $m$ stands for the $n$-dimensional Lebesgue measure.
We set $\mM(\Gamma)=\infty$ when $\mathcal{F}(\Gamma)$ is empty.
The most important property of the modulus is a quasi-invariance;
that is, a  homeomorphism
$f:G\to G'$ between domains in $\Sn$ is $K$-quasiconformal if and only if
\[\mM(\Gamma)/K\le\mM(f(\Gamma))\le K\mM(\Gamma)\]
for  all families of curves $\Gamma$ in $G.$
In particular, $\mM(f(\Gamma))=\mM(\Gamma)$ for a conformal homeomorphism $f.$

For two curve families $\Gamma_1$ and $\Gamma_2$ in ${\Sn}$,
we say that $\Gamma_2$ is minorized by $\Gamma_1$ and
denote $\Gamma_2>\Gamma_1$ if every $\gamma \in \Gamma_2$
has a subcurve which belongs to $\Gamma_1$.
A collection of curve families $\Gamma_j~(j=1,2,\dots)$ is said to be
disjointly supported if there are Borel sets $\Omega_j~(j=1,2,\dots)$
such that all curves in $\Gamma_j$ are contained in $\Omega_j$ and that
$m(\Omega_j\cap\Omega_{j'})=0$ for $j\ne j'.$
Then the following properties of the conformal modulus are fundamental
(see \cite{va} or \cite{gmp}).

\begin{lem}\label{lem:curves}%\cite[Theorem 6.4]{va}
\renewcommand{\labelenumi}{({\arabic{enumi}})}

\begin{enumerate}
\item
If $\Gamma_1 < \Gamma_2$, then $\mM(\Gamma_1) \geq \mM(\Gamma_2)$.
In particular, $\mM(\Gamma_2)\le \mM(\Gamma_1)$ for $\Gamma_2\subset\Gamma_1.$
\item
For a collection of curve families $\Gamma_j~(j=1,2,\dots),$
$$
\mM\left(\bigcup_j \Gamma_j\right)
\le \sum_j\mM(\Gamma_j).
$$
Moreover, equality holds if the collection is disjointly supported.
\end{enumerate}
\end{lem}

A pair $(G,E)$ of a domain $G$ in $\Sn$ and a compact set $E$ in $G$
is called a {\it condenser}.
The {\it capacity of the condenser} $(G,E)$ is
\begin{equation}\label{concapdef}
\capa(G,E)=\mM(\Delta(E,\pa G;G))\,.
\end{equation}
Another equivalent definition makes use of Dirichlet integral minimization
property \cite[Thm 5.2.3]{gmp}.
Here and hereafter,
for sets $E,F,G \subset \Sn$, let $\Delta(E,F;G)$ denote the family of all curves joining
the sets $E$ and $F$ in $G$, and let $\Delta(E,F) = \Delta(E,F; \Sn)$.
Here, a curve $\gamma:[a,b]\to \Sn$ is said to join $E$ and $F$ in $G$
if $\gamma(a)\in E, \gamma(b)\in F$ and if $\gamma((a,b))\subset G.$
For a compact set $E$ in $\Sn$, we write $\capa E =0$  $(\capa E >0)$ if
$\capa(G,E)=0$  $(\capa(G,E)>0)$ for some bounded domain $G$ containing $E$
cf. \cite[7.12]{vubook}.
Note that $\capa(G',E)=0$ for any domain $G'$ containing $E$ if $\capa E=0.$
It is known that $E$ is totally disconnected and has Hausdorff dimension 0 if $\capa E=0\,,$ see
\cite[p.120, Cor.2]{res}, \cite[p. 166, Thm VII.1.15]{ri}.

%{\bf (Reference?)}.

A domain $R$ in $\Sn$ is called a {\it ring} if the complement $\Sn\setminus R$
consists of exactly two connected components, say, $E$ and $F${\black{,}}
and $R$ is often denoted by $R(E,F).$
In particular, $R_{G,n}(s):= R(\Bbn, [s e_1,\infty]),\ s > 1$, is called the {\it Gr\"otzsch ring} and
$R_{T,n}(t):= R([-e_1,0], [t e_1,\infty]),\ t > 0$, is called the {\it Teichm\"uller ring},
where $e_1$ is the unit vector $(1,0,\dots,0)$ in $\Rn.$
%The ring $\mathbb{B}^n \setminus [0, e_1/t], t >1,$ is called the  {\black{\it bounded Gr\"otzsch ring}}.
The capacity of the ring $R(E,F)$ is $\capa R(E,F) = \capa(\Sn\setminus F, E)$
and its modulus is
$$
\mod R(E,F)=\left(\frac{\omega_{n-1}}{\capa R(E,F)}\right)^{1/(n-1)}.
$$
When $R=R(E,F)$ is the standard ring $\{x\in\R^n: a<|x|<b\},$
one has $\mod R=\log(b/a).$
The capacities of $R_{T,n}(t)$ and $R_{G,n}(s)$ are denoted by $\tau_{n}(t)$ and $\gamma_{n}(s)$, respectively.
By \cite[Lemma 5.53]{vubook}, $\tau_{n}: (0,+\infty)\rightarrow (0,+\infty)$
and $\gamma_{n}: (1,+\infty)\rightarrow (0,+\infty)$ are decreasing
homeomorphisms and they satisfy the functional identity
\begin{equation}\label{gamma and tau}
\gamma_n(s)=2^{n-1} \tau_{n}(s^2-1), \quad s >1.
\end{equation}

Here we state a couple of fundamental properties of uniformly perfect sets.
Recall that a ring $R=R(E_1,E_2)$ is said to separate a set $A$ in $\Sn$
if $A\subset E_1\cup E_2$ and $A\cap E_j\ne\emptyset$ for $j=1,2.$
Then the following characterization of uniformly perfect sets is well known
(see, for instance, \cite{aw} for planar case and \cite{gsv} for general case).

\begin{lem}\label{lem:UP}
Let $A$ be a compact set in $\Sn$ with $\card(A)\ge 2.$
Then $A$ is uniformly perfect precisely when there exists a constant $M>0$
such that $\mod R\le M$ for every ring $R$ separating $A.$
\end{lem}

We also note the following simple fact.

\begin{lem}\label{lem:UP2}
Let $G$ be a domain in $\Sn$ for which the complement $C=\Sn\setminus G$
contains at least two points.
Then $\partial G$ is uniformly perfect if and only if so is $C.$
\end{lem}

\begin{pf}
By the previous lemma, it is enough to show that
a ring $R$ separates $C$ if and only if $R$ separates $\pa G.$
Indeed, if a ring $R=R(E_1,E_2)$ separates $C$ then $R\subset G$ and
each $E_j$ meets $C.$
Note that $\Sn\setminus E_2=R\cup E_1$ is a domain.
Choose a point $a$ from $E_1\cap C$ and $z_0$ from $R$
and take a curve $\gamma:[0,1]\to\Sn\setminus E_2$ with
$\gamma(0)=z_0$ and $\gamma(1)=a.$
Then there is a $t\in (0,1]$ such that $\gamma(t)\in\pa G.$
Obviously, $\gamma(t)\in E_1,$ which implies that $E_1\cap\pa G\ne\emptyset.$
Likewise we have $E_2\cap\pa G\ne\emptyset.$
We now conclude that $R$ separates $\pa G.$

Conversely, suppose that a ring $R=R(E_1,E_2)$ separates $\pa G.$
Then $R\subset G$ or $R\cap G=\emptyset.$
If the latter occurs, one component of $\Sn\setminus R,$ say $E_1,$ contains $G.$
Then $E_2\cap\pa G=\emptyset,$ which contradicts the choice of $R.$
Hence the latter case cannot occur.
Therefore, we have shown that $R$ separates $C\,.$ \hfill $\Box$
\end{pf}

For the study of the geometry of the modulus metric below,
we now introduce a new class of conformally invariant domains, M-domains.
The definition of this class  makes use of
the continuum criterion introduced and studied by O. Martio \cite{Mar75}.
The continuum criterion is closely connected with the potential theoretic
boundary regularity of a domain \cite{MarSar77}.

%\begin{defn}% %\end{defn}%
\begin{nonsec}{\bf Definition.}
\label{continuum criterion}
 {\rm
We say that a closed set $C \subset \Rn$ satisfies the {\it continuum criterion}
at $x \in C$ if there exists a continuum
$K \subset \{x\} \cup \left(\Sn \setminus C\right)$ such that
$$\mM(\Delta(K, C;\Sn \setminus C))< \infty.$$
We write $\mM(x,\, C) < \infty$ if this holds,
and otherwise we write $\mM(x,\, C)=\infty.$
 }
\end{nonsec}

We now recall that a continuum is a compact connected set in $\Sn$
containing at least two points.
We note that $\mM(x_0,\,C)=\infty$ if a continuum $C_0 \subset C$ contains $x_0.$
In fact, the sphere $|x-x_0|=r$ meets both $K$ and $C$ for all small enough $r>0$
in this case.
A simple application of the following lemma implies that
$$\mM(\Delta(K,C;\Sn\setminus C))\ge\mM(\Delta(K,C;\Sn))=\infty$$
for every continuum $K$ with $x_0\in K\subset
(\Sn\setminus C)\cup\{x_0\}.$
Here we have used the relation $\Delta(K,C;\Sn\setminus C)<\Delta(K,C;\Sn)$ and
Lemma \ref{lem:curves}.

\begin{lem}[V\"ais\"al\"a $\textrm{\cite[Theorem 10.12]{va}}$]\label{lem:Vai}
Let $0<a<b<+\infty.$
Let  $E$ and $F$ be closed sets in $\Sn$ and suppose that the sphere $|x|=t$ meets both $E$ and $F$
for every $t$ with $a<t<b.$
Then $\mM(\Delta(E,F;\Sn))\ge c_n\log(b/a),$ where $c_n$ is a positive constant depending only on $n.$
\end{lem}

We now define the notion of M-domains.

\begin{defn}\label{def:M-domain}{\rm
A boundary point $x$ of a domain $G \subset \Sn$ is said to
satisfy the {\it M-condition} (relative to $G$) if $\mM (x, \Sn\setminus G) =\infty;$
in other words, the complement $\Sn\setminus G$ does not satisfy the
continuum criterion at $x.$
The domain $G$ is called an {\it M-domain}
if every boundary point $x \in \pa G$ satisfies the M-condition relative to $G.$}
\end{defn}

{\black
By the above observation, a point $x\in \partial G$ satisfies the condition
$\mM(x,\Sn\setminus G)<\infty$ only if  the singleton $\{x\}$ is a connected component of $\partial G.$
On the other hand, any isolated point $x$ of $\partial G$ satisfies
$M(x,\Sn\setminus G)<\infty.$
}

We need the following result in the proof of Theorem \ref{thm:muproper}.
Our proof is similar to that of \cite[Lemma 3.5]{Mar75}.

\begin{lem}\label{lem:ep}
Let $G$ be a domain in $\Sn\,.$
Suppose that a point $x_0\in\pa G\setminus\{\infty\}$ and a continuum $K$ in $G\cup\{x_0\}$
with $x_0\in K$ satisfy the condition $\mM(\Delta(K,\pa G;G))<\infty.$
Then
$$
\lim_{r\to0} \mM(\Delta(K\cap \overline{B}^n(x_0,r),\pa G; G))=0.
$$
\end{lem}

\begin{pf}
If $\pa G=\{x_0\},$ the assertion trivially holds.
Thus we may assume that $\pa G$ contains at least two points.
By the conformal invariance of the capacity, we may assume that $\infty\in\pa G.$
For brevity, we write $\Bber(r)=\overline{B}^n(x_0,r)$ and $S(r)=\pa B(r)$ throughout the proof.
Let $M_0=\mM(\Delta(K,\pa G;G))<\infty$ and choose $r_0>0$
large enough so that $K\subset B(r_0).$
For a decreasing sequence $r_j~(j=0,1,2,\dots)$ with $r_j\to0~(j\to\infty),$
consider the ring $R_j=\{x\in\R^n: r_{j+1}<|x-x_0|<r_j\}.$
We can choose such a sequence so that
$$
c_j:=\capa R_j=\left(\frac{\omega_{n-1}}{\log(r_j/r_{j+1})}\right)^{1/(n-1)}
\quad {\textrm{ satisfies} }  \qquad \sum_{j=0}^\infty c_j<\infty\,.
$$
For instance, for $c_j=2^{-j},$ we define $r_j$ recursively by the formula
$$
r_{j+1}=r_j \exp\left(-\omega_{n-1} \,{c_j^{1-n}}\right)
=r_j\exp\left(-\omega_{n-1}2^{(n-1)j}\right)
$$
for $j=0,1,2,\dots.$
It is obvious that $r_j\to0$ as $j\to\infty$ for this choice.
Let $K_j=K\cap\overline{R_j}$ and denote by $\Delta_j$ the family of curves
joining $K_j$ and $\pa G$ in the set $\{x\in G: r_{j+2}<|x-x_0|<r_{j-1}\}$ for
$j=1,2,\dots.$
Then the families $\Delta_{N+3j}~(j=0,1,2,\dots)$ are disjointly supported
and contained in the family $\Delta(K,\pa G; G)$ for $N=1,2,3,\dots.$
By Lemma \ref{lem:curves} (2) we obtain
$$
\sum_{j=0}^\infty \mM(\Delta_{N+3j})\le\mM(\Delta(K,\pa G; G))=M_0
\quad (N=1,2,3,\dots)
$$
and hence
$$
\sum_{j=1}^\infty \mM(\Delta_{j})\le 3M_0.
$$
For a given number $\eta>0,$ take a large enough integer $N>0$ so that
$$
\sum_{j=N}^\infty \mM(\Delta_{j})<\eta
\quad\textrm{and}\quad
\sum_{j=N-1}^\infty c_j<\eta.
$$
By construction, we easily see that the curve family
$\Delta(K_j,\pa G; G)\setminus\Delta_j$ is minorized by the family
$$\Delta(S(r_j),S(r_{j-1}); R_{j-1})\cup \Delta(S(r_{j+2}),S(r_{j+1}); R_{j+1}).$$
Thus, by Lemma \ref{lem:curves} (1), we obtain
\begin{eqnarray*}
& \ & \mM(\Delta(K_j,\pa G; G)) \\
& \le& \mM(\Delta_j)+\mM(\Delta(K_j,\pa G; G)\setminus\Delta_j)) \\
& \le& \mM(\Delta_j)+\mM(\Delta(S(r_j),S(r_{j-1}); R_{j-1}))
+\mM(\Delta(S(r_{j+2}),S(r_{j+1}); R_{j+1})) \\
&=&\mM(\Delta_j)+\capa R_{j-1}+\capa R_{j+1}.
\end{eqnarray*}
Therefore, we finally have
\begin{eqnarray*}
\mM(\Delta(K\cap\Bber(r_N),\pa G; G))
&\le& \mM(\Delta(\{x_0\},\pa G; G))
+\sum_{j=N}^\infty \Big[\mM(\Delta_j)+c_{j-1}+c_{j+1}\Big] \\
& <& 0+\eta+\eta+\eta=3\eta.
\end{eqnarray*}
Hence we obtain $\mM(\Delta(K\cap \overline{B}^n(x_0,r),\pa G; G))<3\eta$
for $0<r\le r_N.$ \hfill $\Box$
\end{pf}

The next theorem due to Martio \cite[Theorem 3.4]{Mar75} will also
be used in Section \ref{sec:mu-metric}.

\begin{lem}\label{lem:Mar}
Let $G$ be a proper subdomain of $\Sn$ and fix a point $a\in G.$
For a boundary point $x_0$ of $G$ with $x_0\ne\infty,$
set
$$
L(\ep)=\inf_K \mM(\Delta(K,\pa G; G)),
$$
where the infimum is taken over all continua $K$ joining $a$ and the sphere
$S^{n-1}(x_0,\ep)$ in $G.$
Then $\mM(x_0,\Sn\setminus G)=\infty$ if and only if
$L(\ep)\to \infty$ as $\ep\to 0^+.$
\end{lem}

%{\bf This example may be modified into a Cantor-type set.}
It is clear that M-domains are invariant under M\"obius transformations and conformal mappings.
We next give an example of an M-domain which does not have uniformly perfect boundary.

{\black
\begin{nonsec}{\bf Example.} \label{M domain not UP}{\rm
Let $\{s_k\}$ and $\{r_k\}~(k=1,2,3,\dots)$ be two sequences of positive numbers
converging to $0$ monotonically with the following property:
$$
(*) \quad \alpha_k:=s_k-r_k-(s_{k+1}+r_{k+1})>0.
$$
Then the closed balls $\Bber_k=\Bber^n {(s_ke_1,r_k)}$, $k=1,2,\ldots$, are disjoint
because $\dist(\Bber_k,\Bber_{k+1})=\alpha_k>0,$
where $e_1=(1,0,\dots,0)\in\R^n.$
Let $C=\{0\}\cup\bigcup_{k=1}^\infty \Bber_k$ and
$K_0=\{x=(x_1,\dots,x_n)\in\R^n: x_1\le 0\}\cup\{\infty\}.$
Note that the ring $R_k=\{x: r_k<|x-s_ke_1|<r_k'\}$ separates $C,$
where $r_k'=r_k+\min\{\alpha_{k-1},\alpha_k\}.$
Observe that $\alpha_{k-1}\ge\alpha_{k}$ if and only if
$2s_k-s_{k-1}-s_{k+1}\le r_{k+1}-r_{k-1}.$
This condition is fulfilled when $\{s_k\}$ is convex.

(1) The domain $G=\Sn\setminus(K_0\cup C)$ is an M-domain because
every connected component of $K_0\cup C$ is a continuum.
However, $\partial G$ is not uniformly perfect when
$\limsup_{k\to\infty}(r_k'/r_k)=\infty.$
For instance, we can choose a convex sequence $\{s_k\}$ with
$2s_{k+1}\le s_k$ (such as $s_k=2^{-k}$)
and let $r_k=2^{-k}s_k$  for $k\ge 1.$
Then
$$r_{k+1}/r_k=s_{k+1}/(2s_k)\le 1/4, \quad r_k'=2^{k}r_k-(2^{k+1}+1)r_{k+1}$$
and thus
$$
\frac{r_k'}{r_k}\ge 2^k-\frac14(2^{k+1}+1)=2^{k-1}-2^{-2}\to+\infty
$$
as $k\to\infty.$

(2) Let $G=\Sn\setminus C.$
Suppose that the sequence of rings $A_k=\{x: s_k-r_{k}<|x|<s_k+r_k\}$
satisfies the condition $\limsup_{k\to\infty}\mod A_k=\infty.$
For instance, we can take $s_k=2^{-k^2}, r_k=s_k-2s_{k+1}.$
Then $\mM(0,C)=\infty.$
%Indeed, let $K$ be a continuum with $0\in K\subset G\cup\{0\}.$
Indeed, for each $k$ and $t\in (s_k-r_{k}, s_k+r_k),$ the sphere
$|x|=t$ intersects $C$ by definition.
Hence, for any continuum $K$ with $0\in K\subset G\cup\{0\},$
Lemma \ref{lem:Vai} now yields
$$
\mM(\Delta(K, \partial G; G))\ge \mM(\Delta(K,C;\Sn))\ge c_n\,\log\frac{s_k+r_k}{s_k-r_k}
$$
for sufficiently large $k.$
By the assumption, we have $\mM(\Delta(K,\partial G;G))=\infty.$
In this case, the singleton $\{0\}$ is a connected component of $\partial G$
but the condition $\mM(0, \Sn\setminus G)=\infty$ is satisfied.
%$0\in\partial G$ satsifies the M-condition.

(3) Let $G=\Sn\setminus C$ again. Then
$$
\Delta(K_0,C; G)\subset\bigcup_{k=0}^\infty \Delta_k,
$$
where $\Delta_k=\Delta(K_0,\Bber_k;\Sn)$ for $k\ge 1$
and $\Delta_0=\Delta(K_0,\{0\};\Sn).$
Note that $\beta_0:=\mM(\Delta_0)=0.$
Since the ring $R(K_0,\Bber_k)$ contains $R_k$ as a subring, we have
$$
\mM(\Delta_k)=\capa R(K_0,\Bber_k)
\le\capa R_k=\omega_{n-1}(\mod R_k)^{1-n}
=\omega_{n-1}\left(\log\frac{r_k'}{r_k}\right)^{1-n}.
$$

Let $D_k=\{x: |x-s_k|<s_k\}$ for $k\ge 1$ and
$H=\{x: x_1>0\}=\Sn\setminus K_0.$
Then
$$
\mM(\Delta_k)=\capa (H, B_k)\le \capa (D_k,B_k)
=\omega_{n-1}\left(\log\frac{s_k}{r_k}\right)^{1-n}=:\beta_k
$$
for $k\ge 1.$
If $\sum_k \beta_k<+\infty,$ we have
\[
\mM(\Delta(K_0,C; G))\le\sum_{k=0}^\infty \mM(\Delta_k)
%\le \omega_{n-1}\sum_{k=0}^\infty\left(\log c_k\right)^{1-n}
\le\sum_{k=0}^\infty\beta_k<+\infty.
\]
Hence $\mM(0,\partial G)<\infty$ in this case.
For instance, if we choose $s_k$ and $r_k$ so that $r_k=s_k e^{-k^2}$
then $\beta_k=\omega_{n-1}k^{2-2n}$ satisfies the above condition.
Hence, $\mM(0, \Sn\setminus G)<\infty.$
This gives an example of a non-isolated boundary point of a domain
which does not satisfy the M-condition.
}
\end{nonsec}
}

\begin{nonsec}{\bf Open problem.}\label{HdimQ}
{\rm It is well-known that the Hausdorff dimension of the boundary of
a domain with uniformly perfect boundary is positive \cite{JV96}. We do not
know whether the boundary of an $M$-domain has positive Hausdorff
dimension.}
\end{nonsec}

%%%%%%%%%%%%%%%%%%%%%%%%%%%%%%%%%%%%%%%%%%%%%%%%%%%%%%%
%%%%%%%%%%%%%%%%%%%%%%%%%%%%%%%%%%%%%%%%%%%%%%%%%%%%%%%
\section{Modulus metric}\label{sec:mu-metric}
%%%%%%%%%%%%%%%%%%%%%%%%%%%%%%%%%%%%%%%%%%%%%%%%%%%%%%%
%%%%%%%%%%%%%%%%%%%%%%%%%%%%%%%%%%%%%%%%%%%%%%%%%%%%%%%

In this section, we first give a definition of the modulus metric $\mu_G(x,y)$
and its dual quantity $\lambda_G(x,y).$
After that, we will prove Theorems \ref{thm:muproper} and \ref{thm:muminor}.
For further results, we refer to \cite{fe,fe2,Fer97,fmv,LF}, \cite{hkv},\cite{pa1,pa2},
\cite{bpo,ps,z}.

\begin{defn}[$\textrm{\cite[Ch 8]{vubook}}$]\label{modulus metric}{\rm
Let $G$ be a proper subdomain of $\Sn$ and $x, y \in G.$  Then we define
$$\mu_G(x, y) = \inf_{C_{xy}} \mM(\Delta(C_{xy}, \pa G;G)),$$
where the infimum runs over all curves $C_{xy}$ in $G$ joining $x$ and $y$.
We also define
$$\lambda_G(x, y) = \inf_{C_x,C_y} \mM(\Delta(C_{x}, C_y;G)),$$
where the infimum runs over all curves $C_x$ and $C_y$ in $G$
joining $x$ (respectively $y$) and $\pa G.$}
\end{defn}

\medskip

In some special cases, the extremal configurations for the curve
families defining  $\mu_G(x,y)$ and  $\lambda_G(x,y)$ are known.
{\black
Indeed, for the case when $G=\Bn$ and $0\ne x\in\Bn, y=0,$ we have
\begin{equation}\label{eq:Gr}
\mu_{\Bn}(x,0)=\mM(\Delta([0,x], \pa \Bn;\Bn))=\gamma_n(1/|x|)\, ,
%=\frac{2\pi}{\mod(\Btwo\setminus[0,x])}
\end{equation}
and, by the symmetry principle \cite[Thm 4.3.5]{gmp}, with $e=x/|x|,$
\begin{eqnarray}
%\begin{equation}
\label{eq:Tr}
\lambda_{\Bn}(x,0)&=&\mM(\Delta((-e,0], [x,e);\Bn))=
2^{1-n}\mM(\Delta([-\infty,0], [x, e/|x|]; \Sn)) \\ \nonumber
&=&2^{1-n}\mM(\Delta([-e,0], [\frac{|x|^2}{1-|x|^2}e, +\infty]; \Sn))
=2^{1-n}\tau_n(|x|^2/(1-|x|^2))\,,
%\end{equation}
\end{eqnarray}
see \cite[Thm 10.4]{hkv} for details.
Here, we recall that
%$\gamma_n(s)=\capa R_{G,n}(s)$
%=\omega_{n-1}/[\mod R_{G,n}(s)]^{n-1}$
%and $\tau_n(t)=\capa R_{T,n}(t)$ are
the Gr\"otzsch capacity function $\gamma_n(s)$
and the Teichm\"uller capacity function $\tau_n(t)$ are defined by
%where $\Btwo\setminus[0,x]$ is {{ the bounded Gr\"otzsch ring}}.
$$
\gamma_n(s)=\mM(\Delta([0,s e_1], \pa\Bn; \Bn))
\quad\textrm{and}\quad
\tau_n(t)=\mM(\Delta([-e_1,0],[t e_1,\infty]; \Sn))\,,
$$
for $0<s<1$ and $t>0.$

Next we look at the case when $G= \R^n \setminus \{0\}.$
By the definition of $\lambda_G(t e_1, -e_1), t >0,$ there are two natural
choices to connect $t e_1$ and $ -e_1$ with the boundary $ \{0, \infty\}$ of the
domain $G\,,$ either the pair $[t e_1,0), [-e_1,-\infty)$ or  the pair $[t e_1,\infty),
[-e_1,0)\,.$ Therefore
$$\lambda_G(t e_1,-e_1) = \min\{\tau_n(1/t), \tau_n(t) \}$$
and, because $\tau_n: (0,\infty) \to (0,\infty)$ is a strictly decreasing homeomorphism,
for $t>1,$ we have $\tau_n(t)<\tau_n(1)<\tau_n(1/t)$ and thus
$$
\lambda_{\R^n\setminus\{0\}}(t e_1,-e_1)
=\tau_n(t) =\mM(\Delta([-e_1,0),[t e_1,\infty); \R^n\setminus\{0\}))\,, \quad t>1.
$$
%=\frac{2\pi}{\mod(\C\setminus([-1,0]\cup [x,\infty))} = \tau(x)\,.$$
See \cite[p.72]{a} and \cite[pp. 178-181]{hkv} for more details.
}

\medskip

Suppose that $G_1$ and $G_2$ are proper subdomains of $\Sn$ with $G_1 \subset G_2.$
Then for a continuum $C_{xy}$ joining $x$ and $y$ in $G_1$ we have
$\Delta(C_{xy},\pa G_2; G_2)>\Delta(C_{xy},\pa G_1; G_1).$
By Lemma \ref{lem:curves} (1), we further obtain for all $x,y \in G_1$
$$\mu_{G_2}(x,y)\le \mM(\Delta(C_{xy},\pa G_2; G_2))\le\mM(\Delta(C_{xy},\pa G_1; G_1)).$$
Hence $\mu_{G_2} (x, y)\leq \mu_{G_1} (x, y)$.
By definition, the quantities $\mu_G(x,y)$ and $\lambda_G(x,y)$ are both
conformally invariant.
Ferrand \cite{Fer97} proved that $\lambda_G(x,y)^{1/(1-n)}$ is a distance function of $G.$
Thus $\lambda_G(x,y)^{1/(1-n)}$ is often called {\it Ferrand's modulus metric}.
When $n=2$ and $G$ is a simply connected domain in $\Sn$ with $\card(\pa G)\ge 2,$
Ferrand's modulus metric is the same as the modulus metric (up to a constant multiple). Moreover, for $n\ge 2$ there exists  \cite[(9.12), Thm 10.4]{hkv} a constant $c_n >0$ depending only on $n$ such that for all $x,y \in \mathbb{B}^n$
\begin{equation}
\label{murho}
\mu_{\Bn}(x,y) \ge 2^{n-1} c_n \, h_{ \mathbb{B}^n}(x,y)\,.
\end{equation}

\begin{lem}\label{lem:mu-lambda}
Let $G$ be a simply connected hyperbolic domain in $\Stwo=\sphere.$
Then $\mu_G(x,y)=4\lambda_G(x,y)^{-1}.$
\end{lem}

\begin{pf}
{\black
Fix a pair of distinct points $x,y\in G.$
The Riemann mapping theorem asserts that there is a conformal homeomorphism $f:G\to\Btwo
=\{z\in\C: |z|<1\}$
such that $f(x)=0$ and $f(y)=u\in(0,1).$
Since the modulus metric and Ferrand's modulus metric are conformally invariant, we have
$\mu_G(x,y)=\mu_\Btwo(0,u)$
and $\lambda_G(x,y)=\lambda_\Btwo(0,u).$
By \ref{eq:Gr} and \ref{eq:Tr} together with \ref{gamma and tau}, we can write
%It is known (see \cite[Theorem 8.84]{avv} or \cite[Theorem 10.4]{hkv}) that the formulae
$$
\mu_{\Btwo}(0,u)=\gamma_2(1/u)=2\tau_2(u^{-2}-1)
\quad\textrm{and}\quad
\lambda_{\Btwo}(0,u)=\tau_2(1/(u^{-2}-1))/2\,.
$$
In view of the formula $\tau_2(t)\tau_2(1/t)=4$ \cite[5.19 (7)]{avv},
we obtain $\mu_{\Btwo}(0,u)\lambda_{\Btwo}(0,u)=4$ and thus the assertion.
} \hfill  $\Box$
\end{pf}

{\black
We take this opportunity to state the following plausible fact with a short proof.

\begin{lem}\label{lem:top}
Let $G$ be a domain in $\Sn$ such that the complement $F=\Sn\setminus G$
is of positive capacity.
Then there is a positive constant $c(F)$ such that the inequality
\begin{equation}\label{eq:top}
\mu_G(x,y)\ge d_0\min\{q(x,y), c(F)\}
\end{equation}
holds for $x,y\in G,$ where $d_0>0$ is a constant depending only on $n.$
In particular, the modulus metric $\mu_G$ induces the same topology on $G$
as the relative topology on $G$ induced by $\Sn$ with the spherical metric $q.$
\end{lem}

\begin{pf}
The inequality \ref{eq:top} follows from \cite[Theorem 6.1]{vubook}
and implies the inclusion map $(G,\mu_G)\to (\Sn, q)$ is continuous.
In order to show the other inclusion map $(G, q)\to (G,\mu_G)$ is continuous,
we may assume that $G\subset\R^n$ and replace $q$ by the Euclidean metric.
Take an arbitrary point $x\in G$ and choose a small enough number $r>0$
so that $B:=B^n(x,r)\subset G.$
By the domain monotonicity of the modulus metric, we obtain
$$
\mu_G(x,y)\le \mu_B(x,y)=\gamma_n(r/|y-x|), \quad y\in B,
$$
by \ref{eq:Gr}.
Since $\gamma(t)\to 0$ as $t\to+\infty,$ we see that $\mu_G(x,y)\to 0$
as $|y-x|\to 0,$ which proves the required assertion.\hfill  $\Box$
\end{pf}
}
\medskip

We are now in a position to prove the first main result.

\ \newline
\begin{nonsec} {\bf Proof of Theorem \ref{thm:muproper}.}{\rm
The part (i) $\Rightarrow$ (ii) is obvious.
We show now that (ii) implies (iii) by contradiction.
Suppose that $G$ is not an M-domain, namely, $\mM(x_0,\Sn\setminus G)<\infty$
for some $x_0\in\pa G.$
By the conformal invariance, we may assume that $x_0\ne\infty.$
We write $B(r)=B^n(x_0,r)$ and $\Bber(r)=\Bber^n(x_0,r)$ for brevity.
By definition, there is a continuum $K$ with $x_0\in K\subset G\cup\{x_0\}$
such that $M_0:=\mM(\Delta(K,\pa G; G))<\infty.$
Take a point $x_1$ from $K\cap G$ and fix it.
%Let $r_j~(j=1,2,\dots)$ be a decreasing sequence of positive numbers converging to $0$
Let $r_1=|x_1-x_0|$ and $K_1=K.$
For each $x\in K\cap B(r_1)$ and $r\in(0,|x-x_0|),$ let $K_1(x,r)$
be the connected component of $K_1\setminus B(r)$ containing $x.$
Note that $K_1(x,r)$ is a continuum.
By construction, $K_1(x,r)\subset K_1(x,r')$ for $0<r'<r<|x-x_0|.$
We set
$$
C_1=C(x_1,K_1):=\bigcup_{0<r<r_1} K_1(x_1,r).
$$
Then, $C_1$ is connected and,
for $x,y\in C_1,$ we have $x,y\in K_1(x_1,r)$ for some $0<r<r_0.$
%$0<r<\min\{|x-x_0|, |y-x_0|\}.$
In particular, for such a pair of points $x,y$ and $r,$
$$
\mu_G(x,y)\le \mM(\Delta(K_1(x_1,r),\pa G; G))
\le \mM(\Delta(K_1,\pa G; G)).
$$
We also see that $x_0\in\overline{C_1}.$
Indeed, otherwise $\overline{C_1}$ would be a continuum in $K\setminus\Bber(\ep)$
for small enough $\ep>0$ and thus $K_1(x_1,\ep)\supset\overline{C_1}\supset C_1.$
Since $K_1(x_1,\ep)\subset C_1,$ the set $C_1$ would be closed and have a positive distance
to $K\setminus C_1,$ which would violate connectedness of $K.$

Let $K_2$ be the connected component of the compact set $K_1\cap\Bber(r_1/2)$
containing $x_0.$
Since $x_0\in\overline{C_1},$ we have $C_1\cap K_2\ne\emptyset.$
Take a point $x_2$ from $C_1\cap K_2$ and fix it.
As before, set $C_2=C(x_2,K_2).$
Then $C_2\subset C_1\cap K_2.$
Repeating this procedure, we
define sequences of points $x_j,$ continua $K_j$ and connected sets
$C_j$ inductively with the following properties:
\begin{enumerate}
\item
$K_j\subset \Bber(r_1 2^{1-j}),$
\item
$x_j\in C_j\subset C_{j-1}\cap K_j,$
\item
$x_0\in \overline{C_j}\subset K_j,$ and
\item
$\mu_G(x,y)\le \mM(\Delta(K_j,\pa G; G))$ for all $x,y\in C_j.$
\end{enumerate}
In particular, we observe that
$$
\mu_G(x_j,x_k)\le \mM(\Delta(K_j,\pa G; G)),\quad j\le k.
$$
By Lemma \ref{lem:ep}, we have
$$
\mM(\Delta(K_j,\pa G; G))\le \mM(\Delta(K\cap \Bber(r_12^{1-j}),\pa G; G))
\to 0 \quad (j\to\infty).
$$
Hence, we conclude that $\{x_j\}$ is a Cauchy sequence in $(G,\mu_G).$
{\black
Suppose that this sequence is convergent; that is, $\mu_G(x_j, x_\infty)\to 0$ as $j\to\infty$
for some $x_\infty\in G.$
On the other hand, since $|x_j-x_0|\le r_12^{1-j},$ we have $x_j\to x_0$ in $\Sn.$
Lemma \ref{lem:top} now implies that $x_\infty=x_0\in\pa G,$ which is a contradiction.
Therefore, $(G,\mu_G)$ is not complete.
}

Finally, we prove that (iii) implies (i).
If $\capa\pa G=0,$ then $$\mM(\Delta(K,\Sn\setminus G; G))=\mM(\Delta(K,\pa G; G))=0,$$
which is not allowed by condition (iii).
Therefore, $(G,\mu_G)$ is a metric space under the assumption (iii).
Suppose next that the set $X=\{x\in G: \mu_G(x,a)\le r_0\}$ is not compact
for some $a\in G$ and $r_0>0.$
Then there is a point $x_0\in \pa X\cap(\pa G).$
We may assume that $x_0\ne\infty.$
For every $\ep>0,$ there exists a point $x\in X\cap B^n(x_0,\ep).$
By definition of $X,$ $\mM(\Delta(K,\pa G; G))\le r_0$ for a continuum $K$ in $G\cup\{x_0\}$
with $a,x\in K.$
Therefore, under the notation in Lemma \ref{lem:Mar}, we obtain $L(\varepsilon)\le r_0.$
However, the lemma implies that $\mM(x_0,\Sn\setminus G)<\infty.$
By contradiction, we have shown that (iii) implies (i).
\hfill $\Box$}
\end{nonsec}

Next we prove our second result.

%\ \newline

\begin{nonsec} {\bf Proof of Theorem \ref{thm:muminor}.}
{\rm
Since the uniform perfectness is M\"obius invariant (Lemma \ref{lem:UP}),
we may assume that $\infty\in\pa G$ and thus $G\subset\Rn$ and $\diam(\pa G)=+\infty.$

First suppose that the boundary $\partial G$ of $G$ is uniformly perfect.
Lemma \ref{lem:UP2} implies that the complement $E=\Sn\setminus G$ is also uniformly perfect.
By a theorem of J\"arvi and Vuorinen \cite{JV96}, $E$ satisfies the metric thickness condition.
Vuorinen \cite{vu} proved that for such a domain $G$ there exists a constant $b_1>0$
such that for all $x,y\in G$
$$
\mu_G(x,y)\ge b_1 \sj_G(x,y).
$$
Applying \ref{eq:j2j}, we obtain \ref{eq:mu-delta} with $b=b_1/4.$

We next suppose  \ref{eq:mu-delta}.
%that $\mu_G(x,y)\ge b \delta_G(x,y),~ x,y\in G,$ for a constant $b>0.$
Then by Lemma \ref{lem:comparedelta} (3), we have $\mu_G(x,y)\ge b\, j_G(x,y).$
Let $E=\Sn\setminus G$ and
$$
0<c< c_0:= \exp\left[-2\left(\frac{2\omega_{n-1}}{b\log 3}\right)^{1/(n-1)}\right].
$$
We prove now that $\{x: cr\le |x-a|\le r\}\cap E\ne \emptyset$
for every $a\in E\setminus\{\infty\}$ and $r>0.$
Suppose, to the contrary, that $\{x: cr\le |x-a|\le r\}\cap E= \emptyset$
for some $a \in E,~a\ne\infty,$ and $r>0.$
Set $C_1=\{x\in\Rn: |x-a|\le cr\}$ and $C_2=\{x\in\Sn: |x-a|\ge r\}.$
Then the assumption implies that  the set $E$ decomposes into
the two non-empty sets $E_1=E\cap C_1$ and $E_2=E\cap C_2.$
Pick two points $x,y$ from the sphere $S=S^{n-1}(a,\rho)$ so that
$|x-y|=2\rho,$ where $\rho=\sq c \, r.$
We take a curve $C_{xy}^0$ joining $x$ and $y$ in $S.$
Then, by the subadditivity and monotonicity of the modulus
(Lemma \ref{lem:curves}), we obtain
\begin{eqnarray*}
\mu_G(x,y) &\le & {\mM}(\Delta(C_{xy}^0, E)) \\
   &\le& {\mM}(\Delta(C_{xy}^0, E_1))+{\mM}(\Delta(C_{xy}^0, E_2))  \\
&\le& {\mM}(\Delta(S, C_1; G_1))+{\mM}(\Delta(S, C_2;G_2)),
\end{eqnarray*}
where $G_1=\{x: |x-a|<\rho\}$ and $G_2=\{x: |x-a|>\rho\}.$
As is well known \cite[(5.10), (5.14)]{vubook},
\begin{eqnarray*}
{\mM}(\Delta(S, C_1; G_1))={\mM}(\Delta(S, C_2;G_2))
=\frac{\omega_{n-1}}{(\log 1/\sq c)^{n-1}},
\end{eqnarray*}
we have $$\mu_G(x,y)\le \frac{2\omega_{n-1}}{(-\log\sqrt c)^{n-1}},$$
where $\omega_{n-1}$ is the $(n -1)$-dimensional area of $\mb{S}^{n-1}$.
On the other hand, since $d_G(x)\le |x-a|=\rho$ and $d_G(y)\le |y-a|=\rho,$
we obtain $$j_G(x,y)=\log\left(1+\frac{|x-y|}{\min\{d_G(x),d_G(y)\}}\right)
\ge\log\left(1+\frac{2\rho}{\rho}\right)=\log 3.$$
Thus we have $b\log 3\le 2\omega_{n-1}/(-\log\sqrt c)^{n-1},$
that is, $$c\ge \exp[-2(2\omega_{n-1}/b\log 3)^{1/(n-1)}]=c_0,$$
a contradiction.
\hfill $\Box$}
\end{nonsec}

In the case when $G$ is either $\Bn$ of $\Hn,$ the metric $\mu_G(x,y)$ has
the explicit expression in terms of the hyperbolic metric  $h_G$ \cite[Theorem 8.6]{vubook}
\begin{equation} \label{muForm}
\mu_G(x,y)
= 2^{n-1}\,\tau_n \left(\frac1{\sinh^2(\frac12 h_{G}(x,y))}\right)
= \gamma_n \left(\coth^2\Big(\frac{h_{G}(x,y)}2\Big)\right).
\end{equation}

The decreasing homeomorphism $\mu: (0, 1] \rightarrow [0,\infty)$ is defined by
\begin{eqnarray*}\label{muDef}
\mu(r)=\frac{\pi}{2}\frac{{\K}(\sqrt{1-r^2})}{{\K}(r)}\,, \quad {\K}(r)=\int_{0}^{\pi/2} \frac{dt}{\sqrt{1-r^2\sin^2 t}} \,,
\end{eqnarray*}
for $ r \in (0,1)\,, \, \mu(1)=0\,$.
Now the Gr\"otzsch capacity for $n=2$ can be expressed as follows
\begin{equation} \label{GrCap}
\gamma_2(s)= \frac{2 \pi}{ \mu(1/s)}, \quad s >1\,.
\end{equation}
In conjunction with the above relations \ref{muForm}, \ref{GrCap}, when $G$ is the unit disk $\bB^2=\D$ in $\C,$
we obtain the expression
\begin{equation} \label{eq:mu2}
\mu_{\D}(z,w)=
\gamma_2 \left(\frac{1}{\tanh \frac{1}{2}h_{\D}(z,w)}\right)
=\frac{2\pi}{\mu\left({\tanh \frac{1}{2}h_{\D}(z,w)}\right)},
\quad z,w \in \D.\nonumber
\end{equation}

The following estimate will be used later.

\begin{lem}\label{lem:mu}
$$
\mu\left(\tanh x\right)< \frac{\pi^2}{4x},
\quad x>0.
$$
\end{lem}

\begin{pf}
From \cite[(5.29)]{avv}, we note the inequality
$$
\mu (r)<\frac{\pi^2}{4\arth \sqrt[4]{r}}
$$
for $0<r<1.$
Let $v=\left(\tanh x\right)^{1/4}\in(0,1)$ for $x>0.$
Since $0 <\tanh x=v^4<v<1,$ we obtain $x<\arth v.$
Hence,
$$\mu(\tanh x) =\mu \left(v^4\right)< \frac{\pi^2}{4\arth v}<\frac{\pi^2}{4x}.$$ \hfill  $\Box$
\end{pf}

We are now ready to show our third result.

\begin{nonsec}{\bf Proof of Theorem \ref{thm:mumajor}.} {\rm
Assume that $G$ is a M\"obius uniform domain in $\Sn.$
By M\"obius invariance of Definition \ref{MobUnif}, we may assume that $G\subset\R^n.$
By virtue of Lemmas \ref{lem:sigandkap} and \ref{lem:comparedelta},
the uniformity assumption reads
$$
k_G(x,y)\le c \, j_G(x,y),\quad x,y\in G
$$
for a positive constant $c.$
By \cite[Lemma 8.32 (2)]{vubook} (see also \cite[Lemma 10.7]{hkv}) there are
positive constants $b_1, b_2$ depending only on $n$ such that
$$
\mu_G(x,y)\le b_1k_G(x,y)+b_2
$$
for all $x,y\in G.$
In view of Lemma \ref{lem:comparedelta},
we have the required inequality with $d_j=cb_j~(j=1,2).$

Next we assume that the inequality \ref{eq:d} holds for a simply connected
domain $G$ in $\sphere$ with non-degenerate boundary.
We can also assume that $G\subset\C.$
Then, as is well known, the Koebe one-quarter theorem leads to the inequality
$k_G(x,y)\le 2h_G(x,y).$
By the Riemann mapping theorem, there is a conformal homeomorphism $f:G\to\bB^2=\D.$
Since $\mu_G$ and $h_G$ are conformally invariant, we obtain the formula
$$
\mu_G(x,y)=\mu_\D(f(x),f(y))
=\frac{2\pi}{\mu\left({\tanh \frac{1}{2}h_{\D}(f(x),f(y))}\right)}
=\frac{2\pi}{\mu\left({\tanh \frac{1}{2}h_{G}(x,y)}\right)}.
$$
We now apply Lemma \ref{lem:mu} to get
$$
\mu_G(x,y)\ge \frac 4\pi \, h_G(x,y)\ge \frac 2\pi k_G(x,y).
$$
Combining this with \ref{eq:d} and Lemma \ref{lem:comparedelta}, we have
$$
k_G(x,y)\le \frac\pi 2 \mu_G(x,y)\le \frac \pi 2(2d_1 j_G(x,y)+d_2).
$$
Now a result of Gehring and Osgood  \cite{go} implies that $G$ is uniform.
\hfill $\square$}
\end{nonsec}

\begin{nonsec}{\bf Open problem.} \label{HypQ}
{\rm As pointed out above, in the case of planar simply connected domains
the modulus metric can be expressed as a function of the hyperbolic metric.
We do not know, whether for a general hyperbolic planar domain, the hyperbolic metric
has a minorant in terms of the modulus metric.}
\end{nonsec}

%%%%%%%%%%%%%%%%%%%%%%%%%%%%%%%%%%%%%%%%%%%%%%%%%%%%%%%%%%%%%%%%%%%%%%%%%%%%%%%%%
%%%%%%%%%%%%%%%%%%%%%%%%%%%%%%%%%%%%%%%%%%%%%%%%%%%%%%%%%%%%%%%%%%%%%%%%%%%%%%%%%
\section{Application to quasimeromorphic maps}\label{sec:appl}
%%%%%%%%%%%%%%%%%%%%%%%%%%%%%%%%%%%%%%%%%%%%%%%%%%%%%%%%%%%%%%%%%%%%%%%%%%%%%%%%%
%%%%%%%%%%%%%%%%%%%%%%%%%%%%%%%%%%%%%%%%%%%%%%%%%%%%%%%%%%%%%%%%%%%%%%%%%%%%%%%%%
The modulus of a curve family is one of the most important conformal invariants of
geometric function theory which provides a bridge connecting geometry and potential
theory. The modulus is the main tool of the theory of quasiconformal, quasiregular
and quasimeromorphic mappings in $\R^n$ \cite{avv, gmp,va, res, ri,hkv}.
These mappings are the higher dimensional counterparts of the classes of conformal, analytic,
and meromorphic functions of classical function theory, respectively.
We will now apply our results to prove
a M\"obius invariant counterpart of a result of Gehring and Osgood \cite{go} for
quasimeromorphic mappings.

We make use of some basic facts of the theory of quasiconformal, quasiregular,
and quasimeromorphic mappings which are readily available in
\cite{va},  \cite{res}, \cite{ri},  \cite{vubook}.
The first result shows a Lipschitz type property of quasimeromorphic mappings with
respect to the modulus metric. Note that these mappings are locally H\"older-continuous
with respect to the Euclidean metric as some basic examples show \cite[16.2]{va}.

\begin{thm}{\cite[Thm 10.18]{vubook}}\label{thm:qr3} Let $f: G_1 \to G_2$ be a non-constant
$K$-quasimeromorphic mapping where $G_1, G_2 \subset  \overline{\mathbb{R}}^n\,.$
Then for all $x,y \in G_1$, $$\mu_{G_2}(f(x), f(y)) \le K\, \mu_{G_1}(x,y) \,.$$
In particular, $f: (G_1, \mu_{G_1})\to (G_2, \mu_{G_2})$ is Lipschitz continuous.
\end{thm}

%Let  $f:G \to G'\,$ be a homeomorphism between domains $G$ and $ G'$ in $\Sn,~n\ge2\,.$
D.~Betsakos and S.~Pouliasis \cite{bpo} have recently proved that
if $f$ is an isometric homeomorphism between the metric spaces
$$f:(G_1 , \mu_{G_1}) \to (G_2, \mu_{G_2}),$$
then $f$ is quasiconformal and it is conformal if $n=2\,.$
This result gives a solution to a question of J.~Ferrand--G.~J.~Martin--M.~Vuorinen \cite{fmv} when $n=2$.
Very recently this result was strengthened by S.~Pouliasis and A.~Yu.~Solynin \cite{ps} and independently
by X.~Zhang \cite{z}: $\mu$-isometries are conformal in all dimensions $n \ge 2\,.$

%As pointed out in the introduction, it was recently proved that homeomorphisms
%that are $\mu$-isometries are, in fact, conformal mappings.

We next prove a Harnack-type inequality.

\begin{thm} \label{thm:qr4}
Let $f:G_1\to G_2$ be a $K$-quasiregular mapping where $G_1\,, G_2$ are
subdomains of $ {\mathbb{R}}^n\,,n \ge2\,.$
If  the boundary  $\partial G_2$ is uniformly perfect, then the function
\begin{eqnarray*}
u_f(x) :=d_{G_2}(f(x)) = \inf \{|f(x)-z|: z  \in \partial  G_2 \}
\end{eqnarray*}
satisfies the Harnack inequality, i.e. there exists a constant $D_1$ such that
for all $x \in G_1\,,$ and all  $y \in \bar{B}^n(x,d_{G_1}(x)/2)\,,$
\begin{eqnarray*}
\hspace*{50mm} u_f(x) \le D_1\, u_f(y)\,. \hspace*{49mm}(1)
\end{eqnarray*}
Moreover, there exists a constant $D_2$ such that for all $x,y \in G_1$
\begin{eqnarray*}
\hspace*{12mm} k_{G_2}(f(x), f(y)) \le D_2 \, \max \{  k_{G_1}(x,y)^{\alpha} ,  k_{G_1}(x,y)\}\,, \quad \alpha=K^{1/(1-n)}\,. \hspace*{12mm} (2)
\end{eqnarray*}
%Here $k_G$ stands for the quasihyperbolic metric.
\end{thm}

\begin{pf}
Fix $x \in G_1$ and $y \in \bar{B}^n(x,d/2)\,,$ where $d=d_{G_1}(x).$
Then the ring $R=\{z: d/2<|z-x|<d\}$ separates $\{x,y\}$ from $\pa G_1$
and $\mod R=\log 2.$
Therefore, by the definitions of $\mu_{G_1},$
$$
\mu_{G_1}(x,y) \le {{\mM}}(\Delta([x,y], G_1)) \le \capa R
=\omega_{n-1} ({\log 2})^{1/(n-1)}=:M,
%\gamma_n(2) \,.
$$
where we used the relation $\Delta([x,y], G_1)>\Delta(S^{n-1}(x,d/2),S^{n-1}(x,d);R)$
and Lemma \ref{lem:curves} (2).
(A similar estimate is found at \cite[8.8]{vubook}.)
Because $\partial G_2$ is uniformly perfect, it follows from Theorem \ref{thm:muminor}
and Lemma \ref{lem:comparedelta} that
$$
\mu_{G_2}(f(x),f(y)) \ge c \delta_{G_2}(f(x),f(y)) \ge c j_{G_2}(f(x),f(y)) \,.
$$
Next, by Theorem \ref{thm:qr3}
$$
\mu_{G_2}(f(x),f(y)) \le K \, \mu_{G_1}(x,y) \le K \, M\,.
$$
The Harnack inequality (1) with the constant $D_1 =\exp(  K M/2)$ then follows,
because for all $z \in \partial G_2$ \cite[(2.39)]{vubook}
$$
j_{G_2}(f(x),f(y)) \ge \log \frac{|f(x)-z|}{|f(y)-z|}\,.
$$
The proof of (2) follows now from \cite[Theorem 12.5]{vubook}.

\end{pf}

We are next going
to prove the following theorem, which extends a result of F.W. Gehring and B. Osgood \cite[Theorem 3]{go} for quasiconformal mappings. This proof is based on the above Harnack
inequality.

\begin{thm} \label{thm:qr2}
Let $f:G_1\to G_2$ be a $K$-quasimeromorphic mapping where $G_1\,, G_2 \subset \overline{\mathbb{R}}^n\,,n \ge2\,.$ If  the boundary  $\partial G_2$ is uniformly perfect,
then there  exists a constant $d_3>0$ such that for all $x,y \in G_1$
$$
{\sigma}_{G_2}(f(x), f(y)) \le d_3 \, \max \{  {\sigma}_{G_1}(x,y)^{\alpha} ,  {\sigma}_{G_1}(x,y)\}\,, \quad \alpha=K^{1/(1-n)}\,.
$$
\end{thm}

{\black
We prove below in Example \ref{myExamp}  that the uniform perfectness of $G_2$ cannot be dropped from Theorem
\ref{thm:qr2} and the same example also shows that a similar remark
applies to Theorem \ref{thm:qr4}. In this example, the image
domain $G_2$ has one isolated boundary point and cannot therefore
be uniformly perfect.}

\begin{nonsec}{\bf Proof of Theorem \ref{thm:qr2}.}
{\rm Choose M\"obius transformations $f_1, f_2$ such that $0, \infty \in \partial f_1(G_1)$ and
$0, \infty \in \partial f_2(G_2)\,.$ Then
$$g= f_2 \circ f \circ f_1^{-1}: f_1(G_1) \to f_2(G_2) $$
is $K$-quasiregular and by Theorem \ref{thm:qr4} we have
\begin{eqnarray*}
k_{f_2(G_2)}(g(x), g(y)) \le d_3 \, \max \{  k_{f_1(G_1)}(x,y)^{\alpha} ,  k_{G_1}(x,y)\}\,, \quad \alpha=K^{1/(1-n)}\,. \nonumber
\end{eqnarray*}
Because $f_1(G_1), f_2(G_2) \subset \mathbb{R}^n\,,$ we obtain by
Lemma \ref{lem:sigandkap} (3)
a similar inequality for the $\sigma$ metric, with a bit different constants.}
\hfill $\Box$
\end{nonsec}
\bigskip

\begin{nonsec}{\bf Example.} \label{myExamp}{\rm
To show that the condition $\partial G_2$ be uniformly perfect
cannot be dropped from Theorem \ref{thm:qr2}, we consider the analytic function
$g(z)= \exp\left(\frac{z+1}{z-1}\right)$ which maps the unit disk $\mathbb{B}^2$ onto
$\mathbb{B}^2\setminus \{0\}\,.$
Let $G_1=\mathbb{B}^2$ and $G_2=\mathbb{B}^2\setminus \{0\}$,
and let $x_j=(e^j-1)/(e^j+1)$ for $j=1,2,\ldots$.
Then $u_j=g(x_j)=\exp(-e^j).$
The standard formula for the hyperbolic distance \cite[pp.38-40]{be}, \cite[(2.17)]{vubook}
% \cite[11.4]{vubook}
shows that
$$
h_{G_1}(x_j,x_{j+1})=\int_{x_j}^{x_{j+1}}\frac{2dx}{1-x^2}
=2\arth (x_{j+1})-2\arth (x_j)=1
$$
{\black whereas}
$$
k_{G_2}(g(x_j),g(x_{j+1}))=\int_{u_{j+1}}^{u_j}\frac{du}{u}
=e^{j+1}-e^j=(e-1)e^j\rightarrow +\infty
$$
as $j \rightarrow \infty$.
Thus by (1) and (2) of Lemma \ref{lem:sigandkap}, when $j \rightarrow \infty$,
$\sigma_{G_2}(g(x_j),g(x_{j+1}))\rightarrow +\infty$
while $\sigma_{G_1}(x_j,x_{j+1})=h_{G_1}(x_j,x_{j+1})=1$.
This demonstrates that uniform perfectness is needed in Theorem \ref{thm:qr2}.}
\end{nonsec}
%\begin{nonsec}
%\end{nonsec}

%%%%%%%%%%%%%%%%%%%%%%%%%%%%%%%%%%%%%%%
%%%%%%%%%%%%%%%%%%%%%%%%%%%%%%%%%%%%%%%
%%%%%%%%%%%%%%%%%%%%%%%%%%%%%%%%%%%%%%%
\section{Logarithmic M\"obius metric}
%%%%%%%%%%%%%%%%%%%%%%%%%%%%%%%%%%%%%%%
%%%%%%%%%%%%%%%%%%%%%%%%%%%%%%%%%%%%%%%
%%%%%%%%%%%%%%%%%%%%%%%%%%%%%%%%%%%%%%%
\renewcommand{\baselinestretch}{1} \normalsize
In this section we study the logarithmic M\"obius metric
\[\Lam_G(z,w)=\log(1+\delta_G (z,w)) \,, \quad \,z,w \in G\,,\]
%defined by \eqref{logMm}
on a planar domain $G$ in $\sphere = \overline{\R^2}$
and prove Theorem \ref{thm:hypmajor}.
Though the hyperbolic metric $h_G(z,w)$ is majorized by twice
the M\"obius metric $2\delta_G(z,w)$ for an arbitrary hyperbolic domain
$G\subset\sphere$ (see \cite{Seitten99}), the logarithmic M\"obius metric $\Lam_G(z,w)$
is not expected to majorize $h_G(z,w)$ in general.
Indeed, $\delta_G(z,w)$ is Lipschitz equivalent to $h_G(z,w)$ if $\partial G$
is uniformly perfect as we noted in Introduction.
However, the situation is different when $\partial G$ consists of finitely many points.
We now prove the first part of Theorem \ref{thm:hypmajor}.
By using the results from \cite{SZ17} or \cite{SVZ19}, we could obtain
more explicit estimates for the bound $c=c(A).$
However, for brevity, we shall be content with existence of $c>0$ only.

\begin{pf}[Proof of the first part of Theorem \ref{thm:hypmajor}]
Let $A$ be a finite set in $\sphere$ with $\card(A)\ge 3$ and
$G=\sphere\setminus A.$
Since both metrics are M\"obius invariant, we may assume that
$\infty\in A$ so that $G\subset\C.$
We now consider the function
\begin{eqnarray*}\nonumber
\nonumber F(z,w)=\left\{\begin{array}{ll}\nonumber
\displaystyle\frac{h_G(z,w)}{\Lam_G(z,w)} &\quad (z\ne w)\nonumber \\ \nonumber
\vspace{-3mm}\nonumber
&\null \\\nonumber
\displaystyle\frac{\rho_G(z)}{w_G(z)} &\quad (z=w)\end{array}\nonumber
\right.\nonumber
\end{eqnarray*}
on $G\times G.$
Here, $\rho_G(z)$ is the density of the hyperbolic metric on $G$
and $w_G(z)$ is defined in \ref{eq:w}.
Our goal is to find an upper bound of $F(z,w).$
Since the hyperbolic distance is induced by the Riemannian metric $\rho_G(z)|dz|,$
we have
$$
\lim_{w\to z}\frac{h_G(z,w)}{|z-w|}=\rho_G(z)
$$
for $z\in G.$
On the other hand, by definition of the metric $\delta_G(z,w)$
and the property $\log(1+x)=x+O(x^2)~(x\to0),$ we have
\begin{eqnarray}
\lim_{w\to z}\frac{\Lam_G(z,w)}{|z-w|}
&=&\lim_{w\to z}\frac{\delta_G(z,w)}{|z-w|}\nonumber \\
&=&\lim_{w\to z}\frac{m_G(z,w)}{|z-w|}\nonumber \\
&=&w_G(z) \nonumber
\end{eqnarray}
for $z\in G.$
Therefore, we see that the function $F(z,w)$ is continuous on $G\times G.$
Since $\sphere\times\sphere$ is compact, in order to prove that $\sup_{(z,w)\in G\times G}
F(z,w)<+\infty,$ it is enough to prove that
$$
\hat F(\zeta,\omega):=\limsup_{(z,w)\to (\zeta,\omega)}F(z,w)<+\infty
$$
for each $(\zeta,\omega)\in\partial(G\times G).$
Note that $\partial(G\times G)=(\partial G\times G)\cup
(G\times\partial G)\cup(\partial G\times\partial G).$
When $(a,z_0)\in \partial G\times G=A\times G,$
by Lemma \ref{lem:asymp}, we have $\hat F(a,z_0)=1.$
(If $a=\infty,$ with the M\"obius invariance of $F(z,w)$ in mind,
we may consider the inversion $1/z$ to reduce to the finite case.)
Likewise, we can see that $\hat F(z_0,a)=1.$

The remaining case is when $(a,b)\in\partial G\times\partial G.$
We may further assume that $a\ne\infty\ne b.$
If $a\ne b,$ letting $C>|a-b|^2$ be a suitable constant, we have
$$
m_G(z,w)=|a,z,b,w|=\frac{|a-b||z-w|}{|a-z||b-w|}\le \frac C{|a-z||b-w|}
$$
for $z,w$ with $|z-a|<\vep$ and $|w-b|<\vep,$
where $\vep>0$ is a small enough number.
Therefore, taking a fixed point $z_0\in G,$
we have for the same $z,w,$
\begin{eqnarray}
F(z,w)
&\le& \frac{h_G(z,z_0)+h_G(z_0,w)}{\Lam_G(z,w)} \nonumber\\
&\le& \frac{h_G(z,z_0)}{\log\big[1+\log(1+C'/|a-z|)\big]}
+\frac{h_G(z_0,w)}{\log\big[1+\log(1+C'/|b-w|)\big]},\nonumber
\end{eqnarray}
where $C'=C/\vep.$
Taking the upper limit as $z\to a$ and $w\to b,$ with the help of \ref{eq:asymp2},
we finally get $\hat F(a,b)\le 2.$

If $a=b,$ assuming $a=0$ and $\D^*\subset G\subset\C\setminus\{0,1\}$  as before,
we have the estimates
$h_G(z,w)\le h_{\D^*}(z,w)$ and $m_G(z,w)\ge m_{\C\setminus\{0,1\}}(z,w)$
for $z,w\in\D^*.$
Hence, $F(z,w)\le h_{\D^*}(z,w)/\Lam_{\C\setminus\{0,1\}}(z,w).$
The expected claim is now implied by \ref{eq:comb}, which is a consequence of the following lemma.
\end{pf}

Let $E^*:=\{z:\ 0<|z|\leq e^{-1}\}$. For $z_1,z_2 \in E^*$, define
\begin{equation} \label{eq:D}
D(z_1, z_2)=\frac{2\sin(\theta/2)}{\max\{\tau_1,\, \tau_2\}}
+\left|\log\tau_2-\log\tau_1\right|,
\end{equation}
where $\tau_1=\log(1/|z_1|),~ \tau_2=\log(1/|z_2|),~ \theta=|\arg (z_2/z_1)| \in [0,\pi].$
It is known that $D(z_1, z_2)$ is a distance function on $E^*$ (see \cite[Lemma 3.1]{SZ17}).

\begin{lem}
Let $\Omega=\C\setminus\{0,1\}.$
\begin{enumerate}
\item[(i)]
$h_{\D^*}(z_1,z_2)\le (\pi/4)D(z_1,z_2)$ for $z_1,z_2\in E^*.$
\item[(ii)]
$D(z_1, z_2)\leq  M_0{\Lam}_{\Omega}(z_1,z_2)$ for $z_1,z_2\in E^*,$
where $M_0 =2/\log \left( 1+ \log3 \right) = 2.6980\ldots $.
\end{enumerate}
The constants $\pi/4$ and $M_0$ are sharp, respectively.
\end{lem}

\begin{pf}
Part (i) is contained in Theorem 3.2 of \cite{SZ17}.
The sharpness is observed for $z_1=e^{-\tau}, z_2=-e^{-\tau}$ as $\tau\to+\infty.$
We prove only part (ii).
Let $z_1,z_2\in E^*.$
We may assume that $|z_1|\leq |z_2|$ by relabeling if necessary.
Then $|z_j|=e^{-\tau_j}~(j=1,2)$ for some $1 \leq \tau_2 \leq \tau_1 < +\infty.$
We put  $\tau= \tau_2$, $s=\tau_1/ \tau$ and $\vphi=\sin (\theta/2),$ where
$\theta=|\arg(z_2/z_1)|\in[0,\pi].$
Then $s\geq 1$, $0\leq \vphi \leq 1$.
By definition, we have
$$
m_{\Omega} (z_1,z_2) \geq \frac{|z_1-z_2|}{|z_1|}
=\sq{(e^{\tau(s-1)} -1)^2+4\vphi^2 e^{\tau(s-1)}}.
$$
Let $x:=e^{s-1} \geq 1$.
Then
\begin{eqnarray}
{\Lam}_{\Omega} (z_1,z_2)
&\geq&  \log\big[1+\log(1+\sq{(x^{\tau} -1)^2+4\vphi^2 x^{\tau}})\big]
=:{f_1}(\tau,\vphi,x), \quad\textrm{and}\nonumber \\
D (z_1,z_2)&=&\frac{2\vphi}{s\tau}+\log(1+\log x)=:{f_2}(\tau,\vphi,x).\nonumber
\end{eqnarray}
Further let
$$
{f_3}(\tau,\vphi,x):={f_2}(\tau,\vphi,x)-M_0 {f_1}(\tau,\vphi,x).
$$
Then ${f_3}(\tau,\vphi,x)$ is decreasing in $1\le\tau<+\infty$,
and thus ${f_3}(\tau,\vphi,x)\leq {f_3}(1,\vphi,x)$ for $\tau\ge 1.$
By straightforward computations, we have
$$
\frac{\partial^2}{\partial \vphi^2} {f_1}(1,\vphi,x)\leq 0
\quad\textrm{and}\quad
\frac{\partial^2 }{\partial \vphi^2}{f_2}(1,\vphi,x)=0.
$$
Therefore ${f_3}(1,\vphi,x)$ is convex in $0\le \vphi\le 1$.
Since
$$
{f_3}(1,1,x)=\frac{2}{1+\log x}+\log (1+\log x)-M_0 \log(1+\log(x+2)),
$$
it is easy to verify that ${f_3}(1,1,x)$ is decreasing in $1\le x$, which leads to
${f_3}(1,1,x)\leq {f_3}(1,1,1)=0.$
Noting that ${f_3}(1,0,x)=(1-M_0)\log(1+\log x)<0$,
we have ${f_3}(1,\vphi,x)\leq 0$ from convexity,
and thus ${f_3}(\tau,\vphi,x)\leq {f_3}(1,\vphi,x)\leq 0$.
This completes the proof of the required inequality.
To show its sharpness, it is enough to put $z_1=e^{-1}$ and $z_2=-e^{-1}$.
\end{pf}

\begin{rem}{\rm
As an immediate consequence of the lemma, we have the inequality
\begin{equation}\label{eq:comb}
h_{\D^*}(z_1,z_2)\le \frac{\pi}{2\log(1+\log3)}\Lam_{\C\setminus\{0,1\}}(z_1,z_2),\quad
0<|z_1|, |z_2|\le e^{-1}.
\end{equation}
As the reader can observe in the proof, this constant $(\pi/4)M_0\approx 2.11904$
is not sharp.}
\end{rem}

We now complete the proof of Theorem \ref{thm:hypmajor}.

\begin{pf}[Proof of the second part of Theorem \ref{thm:hypmajor}]
Let $G$ be a hyperbolic domain in $\sphere$ with a puncture at the point $a.$
Suppose that $\Phi(\delta_G(z,w))\le h_G(z,w)$ for $z,w\in G.$
By the M\"obius invariance of $\delta_G$ and $h_G,$
we may assume that $a=0$ and that $\D^*\subset G\subset\C.$
Then $m_G(x,-x)\ge |0,x,\infty,-x|=2$ and thus $\delta_G(x,-x)\ge \log 3$
for $0<x<1.$
Therefore, we would have $\Phi(\log 3)\le h_G(x,-x).$
On the other hand, letting $\gamma$ be the upper half of the circle $|z|=x,$
we obtain
$$
h_G(x,-x)\le h_{\D^*}(x,-x)\le\int_\gamma \frac{|dz|}{|z|\log(1/|z|)}
=\frac{\pi}{\log(1/x)}.
$$
Since $\log(1/x)\to+\infty$ as $x\to 0^+,$ we observe that
$h_G(x,-x)\to 0$ as $x\to 0^+,$ which contradicts the above. \hfill
\end{pf}

\section*{Acknowledgments}
The authors are indebted to the anonymous referee for valuable corrections.

%\providecommand{\bysame}{\leavevmode\hbox to3em{\hrulefill}\thinspace}
%\providecommand{\MR}{\relax\ifhmode\unskip\space\fi MR }
% \MRhref is called by the amsart/book/proc definition of \MR.
%\providecommand{\MRhref}[2]{%  \href{http://www.ams.org/mathscinet-getitem?mr=#1}{#2}}
%\providecommand{\href}[2]{#2}

%\bibliographystyle{siamplain}
%\begin{thebibliography}{avv}

\end{document}